%% file: article.tex
\documentclass{article}

\usepackage{algorithm2e}
\usepackage{amsmath}
\usepackage{amssymb}
\usepackage[english]{babel}
\usepackage{color}
\usepackage[T1]{fontenc}
\usepackage{graphicx}
\usepackage[utf8]{inputenc}
\usepackage{subfigure}
\usepackage{url}

\let \epsilon \varepsilon
\newcommand{\rhs}{\mathrm{rhs}\ }
\newcommand{\no}{\nonumber}
\newcommand{\qs}{q^*}
\newcommand{\dotqs}{\dot q^*}

\title{Mode-based derivation of adjoint equations \\- a lazy man's approach}

\author{
Julius Reiss\footnote{julius.reiss@tnt.tu-berlin.de},~~Mathias Lemke\footnote{mathias.lemke@tnt.tu-berlin.de}~~and Jörn Sesterhenn\footnote{joern.sesterhenn@tnt.tu-berlin.de}\\ 
\small{Institut für Strömungsmechanik und Technische Akustik} \\ 
\small{Technische Universität Berlin, Germany}
}

\date{\vspace{-2em}}

\begin{document}

\maketitle

\begin{abstract}
A method to calculate the adjoint solution for a large class of partial differential equations is discussed. 
It differs from the known continuous and discrete adjoint, including automatic differentiation. 
Thus, it represents an alternative, third method. 
It is based on a modal representation of the linearized operator of the governing (primal) system.  
To approximate the operator an extended version of the Arnoldi factorization, the dynamical Arnoldi method (DAM) is introduced.
The DAM allows to derive approximations for operators of non-symmetric coupled equations, which are inaccessible by the classical Arnoldi factorization.  
The approach is applied to the Burgers equation and to the Euler equations on periodic and non-periodic domains.  
Finally, it is tested on an optimization problem.\\[1em]

\end{abstract}

\section{Introduction}
\input{inctex/introduction}

\section{Adjoint Approach \label{secAdjointTheory}}
\input{inctex/theory_adjoint}

\section{The Arnoldi Method \label{sec_arnoldi}}  
\input{./inctex/arnoldi}

\section{The numerical Adjoint \label{secEldidec}} 
\input{./inctex/numericalAdjoint}

\section{Results and Validation \label{secResults} }
\input{./inctex/results}

\section{Conclusion}
\input{./inctex/conclusion}

\subsubsection*{Acknowledgments.}
The authors gratefully acknowledge support by the Deutsche Forschungsgemeinschaft (DFG) as part of collaborative research center SFB 1029 "Substantial efficiency increase in gas turbines through direct use of coupled unsteady combustion and flow dynamics" on project C02. 

\bibliography{local}{}
\bibliographystyle{ieeetr}

\appendix
\input{./inctex/app_calcplans}

\end{document}

%% file: inctex/introduction.tex
Adjoint equations play an important role in numerous areas such as data assimilation, active and passive control, model reduction, sensitivity analysis and optimization.
In particular, the latter is an ubiquitous task in engineering, where often not the system state is sought-after, but instead the factors driving it to a desired goal.
The adjoint yields the dependency of such a goal or objective on influence or design parameters \cite{Jameson1988,GilesPierce2000} and is, thus, a suitable tool for such tasks.
A representative example  is the lift coefficient of a wing (the objective) which could be influenced by a geometry change or by suction of the boundary layer at the wing surface (design parameters).
The full flow field is usually of no interest for designing an airplane.
Rather the most efficient modification of the setup needed to improve the performance towards a desired design goal.
This explains the huge interest of engineers and researchers of various fields in the adjoint method.

The adjoint equations, describing the adjoint field, follow strictly from the equations describing the dynamics or the state of the system (called primal equations throughout the manuscript) and can therefore in principle be unambiguously derived.
Traditionally two approaches are used \cite{GilesPierce2000}: In the \emph{continuous approach} the adjoint equations are derived analytically from the primal equations, and are then discretized to obtain the adjoint solution.
In the \emph{discrete approach} the primal differential equations are discretized firstly. 
Subsequently, the resulting equations (or operators) are transposed, to obtain and solve the adjoint equations. 
This is in particular easy when the operator is explicitly available, as is typical for finite-element computations.  
Alternatively, the discrete adjoint equations are determined by means of automatic differentiation.  
In practice all approaches are successfully used, although in the discrete approach the adjoint solution is discrete-consistent with the system solution, while in the continuous approach inconsistencies of the order of the discretization error are usually present.
For small systems the adjoint equation can also be calculated alternatively by a finite difference approach with reasonable effort \cite{LemkeCaiReissPitschSesterhenn2018}.
If the system state is obtained as an iterative solution the residuum influences the calculation of the adjoint, discussed in \cite{AlbringSagebaumGauger2016}.

The analytic derivation as well as the numerical implementation of the adjoint equations can be challenging. 
Depending on the equations many terms need to be linearized and transformed. 
Furthermore boundary conditions have to be considered carefully.
For practical systems the analytic derivation can be a cumbersome task especially for system-equations describing complex phenomena, e.g.~reactive flows or the coupling of various systems (multi-physics).
Also the complexity of numerical implementations introduced for example by flux limiters or turbulence models can be a challenge.
The discrete automatic differentiation approach promises to circumvent the involved problems of deriving the adjoint equations.
However, in practice problems arise for all terms where a linearization is not straightforward such as switch statements or non-differentiable functions.
Also, a mix of different programming languages can  make the use of such tools prohibitively difficult.

As an alternative to the previously described methods, we propose an on-the-fly construction of the adjoint operator by means of an approximation in a Krylov-like subspace approach.
It builds purely on the evaluation of the primal equations from which an approximative adjoint operator is derived.

The paper is structured as follows:
In Sec.~\ref{secAdjointTheory} the adjoint method is briefly recalled and our notation specified.
In Sec.~\ref{sec_arnoldi} the modal description of operators by the Arnoldi factorization is described and the dynamic Arnoldi method (DAM) is presented.  
This method is used to construct an adjoint operator in Sec.~\ref{secEldidec} including a problem-specific training.  
The full method is applied to test problems in the field of fluid dynamics in Sec.~\ref{secResults}.

%% file: inctex/theory_adjoint.tex
To introduce the adjoint approach a matrix-vector notation is used, analogous to \cite{LemkeReissSesterhenn2014}.
A split in a spatial and temporal part of the governing equations is assumed, which is suitable for many problems of practical interest.
\begin{equation}
  {\cal B}(q) = \partial_t q - \rhs(q) - f= 0 \label{eq_abstract_example_system}
\end{equation}
Wherein the right-hand-side (RHS), contains the spatial terms.
In addition an external forcing $f$ is introduced suitable for representing external forces, e.g.~actuators.
The dependency on the source term $f$ implies an implicit dependency of the system state $q = q(f)$.

The common problem of system parameter variation can also be recast in $f$ as discussed in \cite{LemkeCaiReissPitschSesterhenn2018,GrayLemkeReissPaschereitSesterhennMoeck2017}, but is neglected here for the 
sake of simplicity. 
The adjoint equations of \eqref{eq_abstract_example_system} originates from the scalar-valued objective function
\begin{equation}
  J = J(q(f)),
\end{equation}
which typically encodes a design goal \cite{LemkeReissSesterhenn2014}, but can also be a target for data assimilation \cite{GrayLemkeReissPaschereitSesterhennMoeck2017,Lemke2015}, or can induce a chemical sensitivity analysis \cite{LemkeCaiReissPitschSesterhenn2018}.

If the objective function is non-linearly depending on $f$, either due to a non-linear primal system operator or a direct non-linear dependency of $J$ on $q$, a linearized formulation is needed.
The linearization of the primal equation (\ref{eq_abstract_example_system}) around the base state $q \to  q + \delta q$  and $f \to  f + \delta f $ can be written as  
\begin{align}
  \dfrac{\partial {\cal B}}{\partial q} \delta q + \dfrac{\partial {\cal B}}{\partial f} \delta f =  
  \underbrace{ \left( \partial_t  - \frac {\partial \rhs(q)}{\partial q} \right)}_{=B } \delta q - \delta f =0  . \label{defB}
\end{align}
The linearization of $J$ can be written as a scalar product 
\begin{equation}
  \delta J = \underbrace{\left( q - q^{\mathrm{(mod)}}\right)^T}_{= g^T} \delta q. \label{eq_abstract_example_linearised_objective}
\end{equation}
The vector space of $\delta q$ and by this the corresponding scalar product spans space \emph{and} time. 

We seek to calculate the change of the objective function $\delta J$ by a change of $f$, under the constraint \eqref{defB}. 
The gradient is derived by the Lagrangian formalism by adding \eqref{defB} to \eqref{eq_abstract_example_linearised_objective} using a multiplier $q^*$
\begin{eqnarray}
  \delta J &=& g^T \delta q - {q^*}^T \underbrace{\left( B \delta q +  \delta f \right)}_{=0} \label{eq_lagrange_ansatz_1_line}\\
	   &=& \delta q^T \left( g - B^T q^* \right) - {q^*}^T  \delta f.
\end{eqnarray}
Since the value of \emph{adjoint variable} $q^{*T}$ is arbitrary, one can determine its value such that the change in the objective function becomes independent of $\delta q$. 
This is the central point of the adjoint approach, since by this the typically computationally expensive solution of \eqref{defB} for every change $\delta f$ can be avoided by demanding that the adjoint state $q^*$ is the solution of the, so-called, adjoint equation
\begin{equation}
  g - B^T q^* = 0, \label{eq_adjoint_equation}
\end{equation}
one obtains 
\begin{equation}
  \delta J = - {q^*}^T  \delta f,
\end{equation}
which is independent of $\delta q$.
 
We apply this derivation to our problem by  defining the linear spatial operator  
\begin{align}
A = {\frac {\partial \rhs(q)}{\partial q}}.
\end{align}
as the spatial part of $B$ in \eqref{defB}.
To derive the adjoint equation for our form, the linearized equation is added to the linearization of the objective function
as before 
\begin{equation}
  \delta J =  \int \limits_{t_0}^{t_{\mathrm{end}}} g^T \delta q ~\mathrm d t
  \int \limits_{t_0}^{t_{\mathrm{end}}} {q^*}^T \underbrace {\left( \partial_t \delta q - A \delta q -  \delta f \right)}_{=0} \mathrm d t. \label{eq_lagrange_ansatz}
\end{equation}
Note the minor change in the definition of the scalar product which covers only the discrete spatial part, since the time is explicitly treated.  
The resulting equation \eqref{eq_lagrange_ansatz} is rearranged by means of an integration by parts of the term
\begin{equation}
  \int \limits_{t_0}^{t_{\mathrm{end}}} {q^*}^T \partial_t \delta q \mathrm d t = \left[ {\delta q}^T q^* \right]_{t=t_0}^{t = t_{\mathrm{end}}} - \int \limits_{t_0}^{t_{\mathrm{end}}} \delta q^T \partial_t q^* \mathrm d t 
\end{equation}
and, in particular, transpose of the operator $A$
\begin{equation}
   \delta J = \int \limits_{t_0}^{t_{\mathrm{end}}} \delta q^T \left( g + \partial_t q^* + {A}^Tq^* \right) \mathrm d t - \left[ {\delta q}^T q^* \right]_{t=t_0}^{t = t_{\mathrm{end}}} + \int \limits_{t_0}^{t_{\mathrm{end}}} {q^*}^T \delta f  ~\mathrm d t.
\end{equation}
The dependency of $\delta q$ is removed by demanding 
\begin{equation}
  g + \partial_t q^* + {A}^Tq^* = 0, 
\end{equation}
resulting in the adjoint equation
\begin{equation}
  \partial_t q^* =  -{A}^Tq^* - g. \label{eq_adjoint_system}
\end{equation}
As the temporal boundary term (initial condition) also depends on $\delta q$, it needs to vanish as well.
\begin{equation}
  - \left[ {\delta q}^T q^* \right]_{t=t_0}^{t = t_{\mathrm{end}}} = - \left[ {\delta q}^T q^* \right]_{t = t_{\mathrm{end}}} + \left[ {\delta q}^T q^* \right]_{t=t_0} = 0
\end{equation}
The right term vanishes as the initial condition for $q$ is fixed and $\delta q(t_0) = 0$ holds. 
To cancel the remaining term, the adjoint state is chosen as $q^* (t_{\mathrm{end}}) = 0$, as the primal state $\delta q (t_{\mathrm{end}})$ is arbitrary at the final time.
In other words, the temporal (initial) condition of the adjoint system is given at final time.
In general, the adjoint system is well-posed only if the adjoint initial state is defined at the end of the computational time and the system is integrated backwards in time, see \cite{GilesPierce2000}.

The expression ${q^*} $ can be interpreted as the sensitivity of the objective function with respect to source $f$, since
\begin{equation}
  \dfrac{\delta J}{\delta f} = {q^*}.  \label{eq_adjoint_sensitivity}
\end{equation}

%% file: inctex/arnoldi.tex
As seen in the previous section, the derivation of the adjoint system would be trivial if the Matrix $A$ was explicitly given, since deriving the adjoint reduces in this case to a simple transpose.
Unfortunately, for large systems, as for flow simulations, this matrix is usually not explicitly available.  
Only the \emph{evaluation of a RHS of a given (spatial) state} is available.  
Methods which rely on the evaluation of a RHS instead of needing the explicit matrix are called matrix-free methods. 
They often build on approximating the matrices by a few, representative modes.

If one chooses  $m$  orthogonal vectors of size $n$ forming the orthogonal $(n,m)$-matrix $V$ the action of a $(n,n)$-matrix $A$ can be described by
\begin{equation}
  V H + r  = A V,
\end{equation}
where $H$ is a $(m,m)$ matrix and $r$ is a residuum orthogonal to $V$. 
By neglecting the residuum the approximation of $A$ results in 
\begin{equation}
  V H V^T = \tilde A. \label{tildeH}
\end{equation} 
The vectors forming $V$ are often referred to as modes. 
In many practical problems it is possible to work directly with a much smaller matrix $H$ instead of $\tilde A$, e.g.~\cite{SchulzeSchmidSesterhenn2009}. 
However, if the matrix $A$ has (near) full rank it is not possible to find  a few modes $V$, so that (\ref{tildeH}) is in general a good approximation.
In this case, the choice of the vectors $V$ is crucial and problem dependent. 
Often the Krylov space $\mathcal{K}_m$ spanned by 
\begin{equation}
  \mathcal{K}_m = span\{ v , Av, A A v , \dots, A^{m-1} v  \} 
\end{equation} 
with $v$ as a suitable initial vector. 
The space consists of applying the operator $A$ in $m-1$ subsequent steps. 
The calculation of the Krylov base and the implied matrix $H$ is done by the Arnoldi factorization \cite{Saad2000,Sorensen2002}, which is a matrix-free method.
If more than one Krylov subspace is to be approximated the block-Arnoldi method calculates the space 
\begin{equation}
  \mathcal{K}_m = span\{ v_k , Av_k, A A v_k , \dots, A^{m-1} v_k  \} 
\end{equation}  
where $v_k$ is a matrix formed by $m$ vectors of length $n$. 

Note, that in both Arnoldi methods the initial vector $v$ or set of initial vectors $v_k$ determines $V$ and by  this the quality of the approximation. 
There are only limited  possibilities to influence the Krylov space, such as restarting with a vector $\bar v$ typically formed from $V$.  

\subsection{Dynamic Arnoldi Method}
\label{sec_dam} 
We find further down that neither the Arnoldi nor the block-Arnoldi method is flexible enough for the application in mind.  
To allow a direct intervention we now define the so termed \emph{Dynamic Arnoldi Method} (DAM).  
It allows to choose in each step freely new vectors to expand the mode-set $V$.  
These vectors are chosen from some initial set or taken from previous calculated application of $A$ and can also be modified before the application of $A$.  
This choice is governed by a set of predefined rules termed the \emph{calculation  plan}. 
By selecting suitable {calculation plans} one recovers the classical Arnoldi or the block-Arnoldi method.

At the core of the DAM is an update step of the relation 
\begin{equation}
 P^m +  V^m \bar H^m = A V^m.
\end{equation}  
The operator $A$ has dimension $(n,n)$  and  the  matrices $P^m, V^m$  have dimension $(n,m)$ and   $\bar H^m$ $(m,m)$.  
The matrix $\bar H $ is \emph{not} the same as $ H$, introduced above, the connection is provided further down. 
The $(n,m)$ matrix $P$ was introduced  acting as a pile for results which, at the time of update, cannot be described by $V$. 
It extends the residuum $r$ defined above. 
The update is done by adding a vector $v$ to $V$ and expanding $P$ and $\bar H$:\\

\paragraph{DAM Update Routine}~\newline
\begin{algorithm}[H]
 	\KwData{ $q^{m+1} ,  V^m,\, \bar H^m, \,  P^m ;~A$ }
 	\KwResult{ $V^{m+1}$, $\bar H^{m+1}$, $P^{m+1}$    }
 	\# orthogonalize input \;   
 	$\alpha = (V^{m})^T\cdot q^{m+1} $ \;
 	$ v^{m+1} = q^{m+1} - (V^{m})\cdot\alpha $\; 
 		\eIf{ $|v^{m+1}|_2 > \epsilon $  }{
 			\# input has linear independent part, append to $V$, apply $A$\;
 			$V^{m+1} =  \mathrm{appendColumn} (V^m, v^{m+1} )$\; 
 			$ w = A \cdot v^{m+1} $  \;
 			$ \beta = (V^{m+1})^T\cdot w  $ \;
 			$ w = w - (V^{m+1})\cdot\beta $\; 
 		}{
 		\# input is linear dependent, append zero column \;  
 	    $V^{m+1} =  \mathrm{appendColumn} (V^m, 0 )$\; 
 		 			
 		$w       = 0  $ \# no need to do matrix multiplication\; 
  	}
  	\vspace{1em}
   $\bar H ^{m+1} =  \mathrm{appendColumn} ( \bar H^m, \beta )$\; 
   $ P  ^{m+1} =  \mathrm{appendColumn} ( P^m, w )$\; 
  ~\\    	
 \caption{The update routine part of the dynamic Arnoldi.}
\end{algorithm}~\\
This algorithm allows to prescribe an arbitrary sequence of $q^{m}$ to multiply with the matrix $A$ and to store the result in $V^m$, the pile $P^m$ and $\bar H$. 
The vector $ q^{m+1}$ is entitled \emph{test vector} in the following.  

\noindent Some remarks are worthwhile:

The operator $A$ is later defined purely in terms of applications of the RHS on a vector.
The method is a matrix-free method building on the application $A\cdot q$ only.    
 
The pile $P$ is in general not orthogonal to $V$.
If, for example, a vector from  $P$ is used as a new test vector, it is added to $V$, but no action is taken to remove it from the pile $P$ in the update step. 
It could be included in each step, but for simplicity we remove the part of $P$ resolved by $V$ (if necessary) only in a final step of the main routine, see below.

Further, zero modes are added to $V$, which seems unnecessary, since it does not enlarge the space spanned by $V$. 
However, dropping those entries would changed the position of results in $P$, which  would complicate the calculation plan discussed below. 

The \emph{DAM Update Routine} is in the heart of the dynamic Arnoldi method which prescribes a sequence of test vectors $q^m$ based on the input vectors and the pile vectors. 
A set of  $l$  input vectors is  provided as an $(n,l)$ Matrix $U$.
To control the DAM  we define the  calculation plan  $\mathcal{C} $.  
It consists of lines $l$
\begin{equation}
  \mathcal{C}^l = ( source ,\; index ,\; mod ).
\end{equation} 
In this work the field $source$ contains the value 'I' for taking the next test vector $q^{m+1}$ from the input $U$ and 'P' for taking it from the pile $P$. 
The $index$  is simply the number of the vector within the input or the pile. 
The modification identifier $mod$ is specific for the  application of this report. 
It is a mask which  allows to shuffle  the different fields within a test vector.  
E.g. a compressible flow in one dimension  may be described by the three fields of density, velocity and  pressure $(\rho,\, u,\,  p)$, 
a modifier  entry $M: (\rho,\, u,\,  p) { \overset {(2,1,3)} \longrightarrow }  (u, \rho , p ) $ would exchange the first and the second field, i.e. the density and the velocity. 
This will prove useful later.
An entry '0' will simply set the field to zero. As an example $(\rho,\, u,\,  p) { \overset {(2,0,0)} \longrightarrow }  (u, 0 , 0  )$ would write the second field to the first an set the others zero.   
These modifications are discussed further  in Sec.~\ref{ResultEuler}.  
The DAM  is given by

\paragraph{Dynamic Arnoldi Method}~\newline
\begin{algorithm}[H]
\KwData{ $ U ;\, \mathcal{C}, \, A$ }
\KwResult{ $V^{m}$, $\bar H^{m}$, $P^{m}$    }
\# initialize 
$V^{0},\, \bar H^{0}, \,  P^{0}$ = empty  \;  
\For{ $C^l = l^{th}$ line in  all lines of $\mathcal{C}$ }{
	\# unpack information \;  
$	source, index, mod = l_k$ \;  
    \# choose source \; 
    \Switch{source}{
    \Case{I}{
    	 $q = U^{index}$ \;  
    	}
    \Case{P}{
    	$q = P^{index}$ \;  
    	} 	 
    }
    \# modify as described in text  \; 
$    q  { \overset {mod} \longrightarrow } q^m  $\; 
    \# update by  Algorithm 1
    $ [ V^{m}, H^{m} , P^{m} ]= DAM\_Update(q^m ,  V^{m-1},\, \bar H^{m-1}, \,  P^{m-1}  ; A  )  $ 
	}
	\# calculate $H$ matrix by adding parts of $P$ described by  $V$  
	$H^m = \bar H^m + (V^m)^T P^m  $ ;\\[1em]
\caption{ The dynamic Arnoldi Method }
\end{algorithm}~\\

For a given calculation plan this allows to create a specific approximation of the matrix $A$. 
The standard Arnoldi is obtained by simply using the last result of the pile  as the next test vector $q^{m+1}$. 
The block-Arnoldi is recovered by  providing a set of input vectors in sequence,  followed by using the resulting pile vectors repetitively as the next test vectors.

%% file: inctex/numericalAdjoint.tex
We aim to approximate the application of the transpose of $A$ to the current state of the adjoint variable $q^*$
\begin{equation}
  \dotqs = A^T q^*
\end{equation}  
by evaluation of $A$ on test vectors. 
To investigate the needed test vectors assume that not only $\qs$ but also the result of the adjoint right-hand-side $\dotqs = A^T \qs $ is known. 
Creating $ V   = (\qs, (\dotqs)^\perp)$, where $ (\dotqs)^\perp$ is the part of $\dotqs$ orthogonal to $\qs$, we obtain the approximation for $A$ 
\begin{equation}
  VH = AV + r , 
\end{equation}
where the residuum $r$ is orthogonal to $V$. 
From this, the approximation for action of the adjoint operator is  $\overset{\sim}{\dotqs} = VH^TV^T \qs$.
It is easy to see that this approximation is exact.    
Completing $V$ to a complete base does not change the approximation since all added base vectors are orthogonal to $ \qs $ and $\dot\qs$. 
Thus, the matrix elements of $H$ contain already all information for the connection of $\qs$ and $\dotqs$.      
A sufficient criterion for a good approximation is how well $\qs$ and $\dot \qs$ are included in $V$.   
It is therefore our goal to utilize the dynamic Arnoldi method to have $\dotqs$ in $V$. 
It is trivial to have $\qs$ by providing it as a first test vector in the dynamic Arnoldi. 
With this choice  the approximation of $\dotqs$ is 
\begin{equation}
  \dot\qs = VH^T V \qs = Vh_1  \label{adjointRhsDAM}
\end{equation}
where $h_1$ is the first line of H. 
Note, that the first line of $H$ is the same as the first line of $\bar H$ so that $H$ does not need to be calculated (albeit cheap).
To have a good approximation of $\dot\qs$ is the guiding principle to determine a calculation plan.  
The DAM is provided with the input of the adjoint solution, the primal solution and possible the previous calculated right-hand-side $U=(\qs,q,(\dotqs)^{\mathrm{old}})$.   
  
\subsection{Linearization of the RHS \label{subsecLinOfRHS}}  

All equations considered further down  are non-linear. 
The linearization is derived as a Frech\'{e} derivative, 
\begin{equation}
  A \cdot q \approx \dfrac{\rhs \left(q_0 + \epsilon q  \right) - \rhs \left(q_0  \right)}{\epsilon}. \label{eq_frechet_of_f}
\end{equation} 
The choice of $\epsilon$ influences the approximation quality, here, we chose the parameter $\epsilon = \sqrt{\epsilon_{\mathrm{machine}}} $, where $\epsilon_{\mathrm{machine}} $ is the machine precision of the used computer.  

This formula can in principle be also used to calculate all elements of the desired matrix $A$ by setting $q$ to all unit vectors of the discrete space dimension \cite{LemkeCaiReissPitschSesterhenn2018}.
However, this becomes prohibitive expensive for large systems. 
It is, however, used in the numerical examples to obtain a reference solution. 
It is also used in the current version of the training of the mode-based adjoint method. 
We referrer to it as expensive method.

\subsection{Training of the method \label{sec_training_of_the_method}}
The calculation plan is created by a training. 
It aims at approximating the adjoint operator $A^T$ with a  minimal number of evaluations of the primal RHS.  
For this, the result of the adjoint RHS (\ref{adjointRhsDAM}) for different training-plan modifications is compared with a reference result. 
The later is, for now, created by the expensive method described in the previous section.
The comparison is done for a representative set of $\qs$, since by using only one, the calculation plan might reflect a specific property of this state not valid in general. 
Here, we use a set of simulation snapshots created with the expensive method. 
In detail, the training plan is initiated with one line applying $A$ to an unmodified $\qs$, which is a necessary mode as discussed at the start of section \ref{secEldidec}. 
Then, in a greedy approach, all possible next training plan lines are tried and the one which reduced the combined error 
for the whole $\qs$-set, is appended. 
The last step is repeated until the error is below a prescribed value. 
Different resulting calculation plans for the test cases are given in appendix \ref{secCalcPlan}.  

The order of steps prescribed by the calculation plan has an impact on the error reduction of a given step, since parts which are represented by previous test vectors does not contribute to this error reduction anymore. 
Also the entries in $P$ change since the new entry is orthogonal to $V$ at the step of calculation. 
This introduces a strong non-linearity by which the greedy strategy can become suboptimal. 
Indeed, we observed for simple cases, where the calculation plan can be constructed by hand, that the training delivers an inferior solution, i.e.~one with more than minimal evaluations of the 
RHS\footnote{A simple example is given by a case where one mode approximates  $\dotqs$ very good and two other modes less good but combine to a perfect representation.
The strategy will pick first the first mode and than the next two which would suffice alone.}.
To reduce this non-linearity it was chosen not to orthogonalize the pile vectors (by introducing $P=\bar P R$ ) or remove the parts contained in $V$, since the pile vectors are less dependent on the history of the updates. 
We observed, thereby, a much more robust (albeit not optimal) training outcome. 
Other strategies like genetic algorithms  or Monte Carlo tree search \cite{BrownePowleyWhitehouseLucasCowlingRohlfshagenTavenerPerezSamothrakisColton2012} are likely to improve this but are out of scope of this paper.      

The reference used in the training is derived in an expensive manner. 
This might still be possible for practical problems since the training can be done for very small systems still capturing the system dynamics and yield a training plan which is assumed to work independent of the discretization size.
In the end, the training plan should reflect the mathematical  structure of the problem. 

If the smallest system is too big for this approach one could try, in a Monte Carlo fashion,  small disturbances at a set of random locations and integrate the original equation forward in time. 
The influence of this should  correctly be predicted by the numerical adjoint replacing the full reference solution in the training. 
This, again, is out of scope of this report, dedicated to the principle idea of the numerical adjoint.

%% file: inctex/results.tex
In the following we present numerical results of the proposed method.
For this, the numerical adjoint of the non-linear Burgers-equation (B) and the Euler-equations (E) is analyzed in different setups.
The Burgers-equation allows to investigate  effects of the non-linearity of the primal equation on the adjoint as well as effects of a friction term, 
which is a symmetric operator in contrast to the near skew-symmetric transport term.
The Euler-equations are employed to discuss additional difficulties of non-symmetric coupled equations as well as a non-periodic spatial discretization.

All results are compared with a discrete adjoint obtained by calculating the linearized operator by a finite difference method explained  in Sec. \ref{subsecLinOfRHS}, 
see also \cite{LemkeCaiReissPitschSesterhenn2018}. 
This gives an reference solution without additional discretization discrepancies, but it is prohibitive expensive for other than small test systems.

All results  shown are normalized in terms of space and time.
   
\subsection{\textsc{Burgers} Equation}
The Burgers equation is given by 
\begin{equation}
  \partial_t u  + \partial_x \left( \frac{u^2}{2} \right) = \mu \partial_x^2 u, \label{eq_burgers}
\end{equation} 
with a scalar transported  quantity $u$ and a friction constant $\mu$, which is set to zero for the friction-less cases. 
The equation is spatially discretized  by  central finite differences.  
A standard central fourth order derivative is used, also for the second derivative by applying the first derivative twice. 
The periodic computational area of length $2\pi$ is resolved by 128 equidistantly distributed points, if not stated otherwise.
The time integration is realized by a standard Runge-Kutta scheme of fourth order.
A total number of $256$ time steps is used to resolve the time span from $t_0$ to $t_{\mathrm{end}}$, which corresponds to about one convectional length for all setups using a CFL condition of $0.5$.
The numerical, Arnoldi-based adjoint and the reference adjoint system are discretized in the same manner.

To aid the discussion of the results, the corresponding adjoint of \eqref{eq_burgers} is analytically derived as
\begin{equation}
  \partial_t u^* + u_0 \partial_x u^* = -\mu \partial_x^2 u^* \label{eq_burgersAdjoint}
\end{equation} 
with $u^*$ as adjoint variable.
All adjoint computations, either based on the proposed method or on an analytical derivation, are initialized by means of a Gaussian disturbance of form
\begin{equation}
 u^*(x,t_{\mathrm{end}}) = \frac{1}{2} \cdot \exp\left(-\frac{(x - x_0)^2}{(15 \Delta x)^2} \right),
\end{equation}
with $x_0$ as center of the computational domain and the grid spacing $\Delta x$. 
This initial condition is set at the end of the computational time $t_{\mathrm{end}}$, as the adjoint is integrated backwards in time.
This mimics an action of a non-zero source term $g$ in the adjoint equation at the last time step ($t_{\mathrm{end}}$).    

The  calculation plans for the different setups, required for the mode-based adjoint method, can be found in Tbl.~\ref{app_tbl_calcplans_burgers} found in the appendix.

\begin{figure}
	\centering
	\includegraphics[width = .49\textwidth]{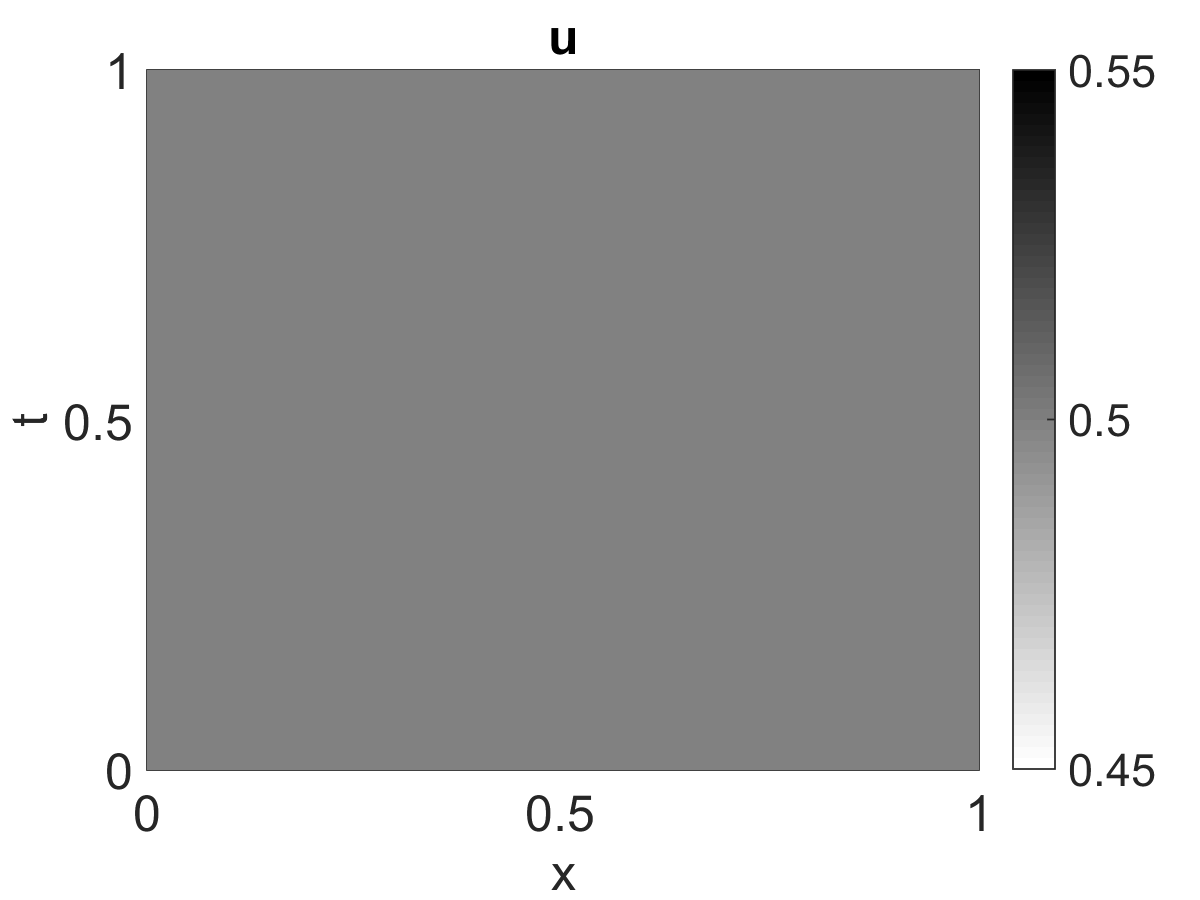}
	\includegraphics[width = .49\textwidth]{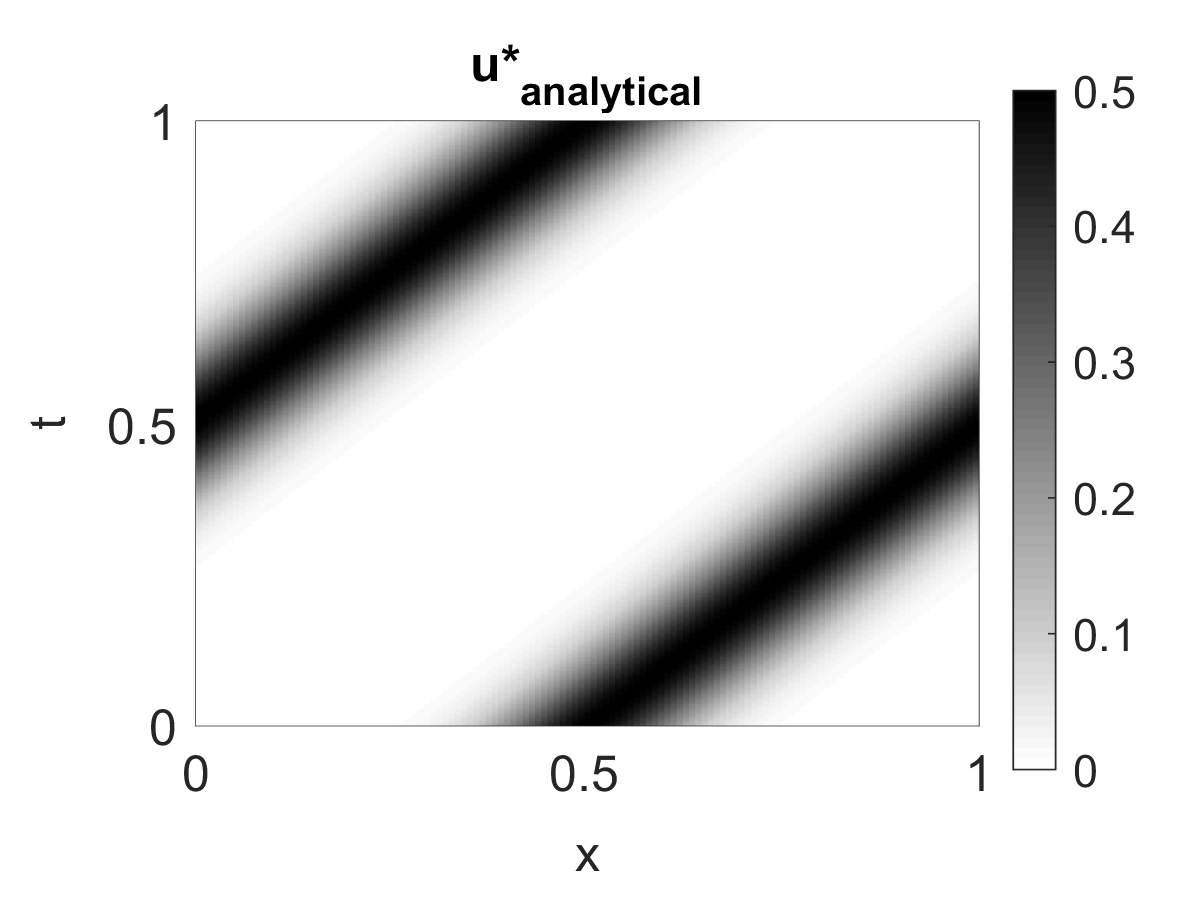}\\
	\includegraphics[width = .49\textwidth]{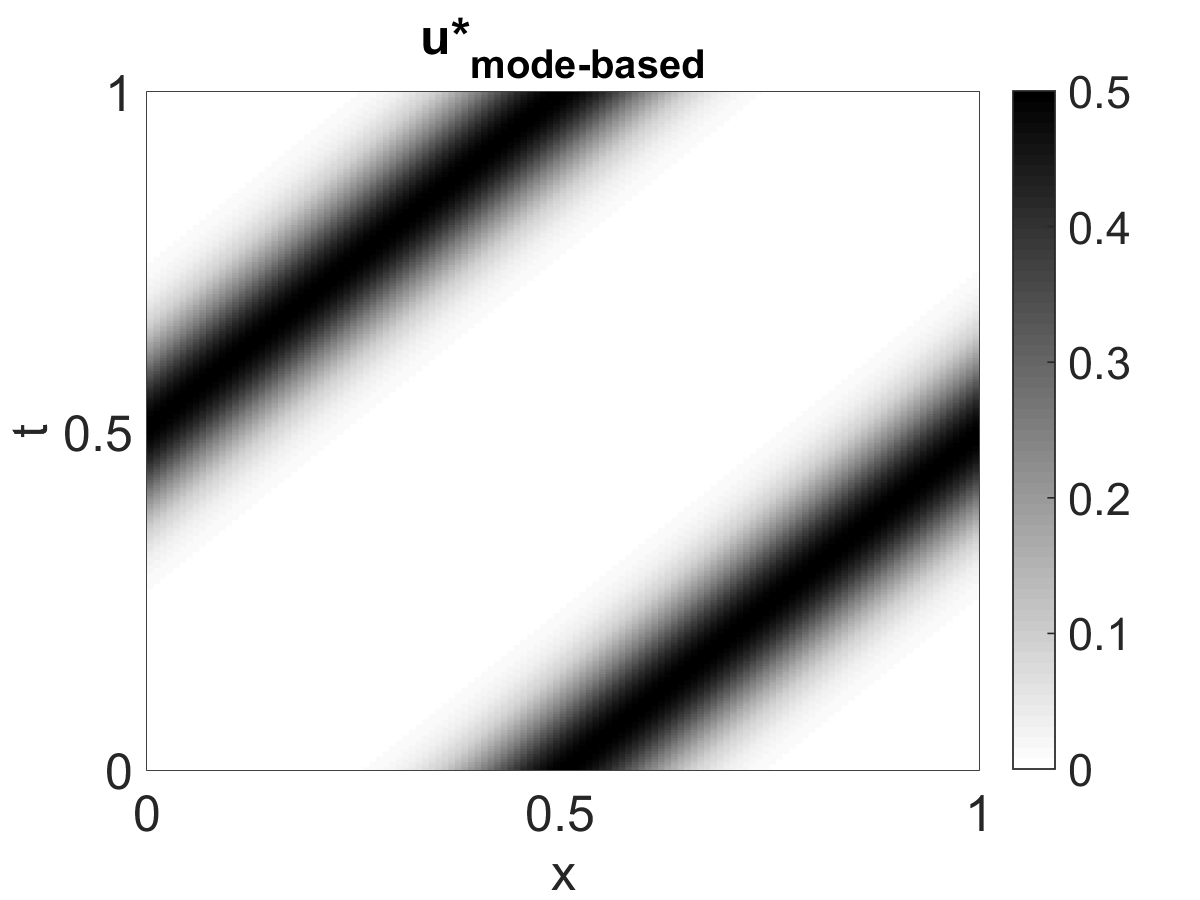}
	\includegraphics[width = .49\textwidth]{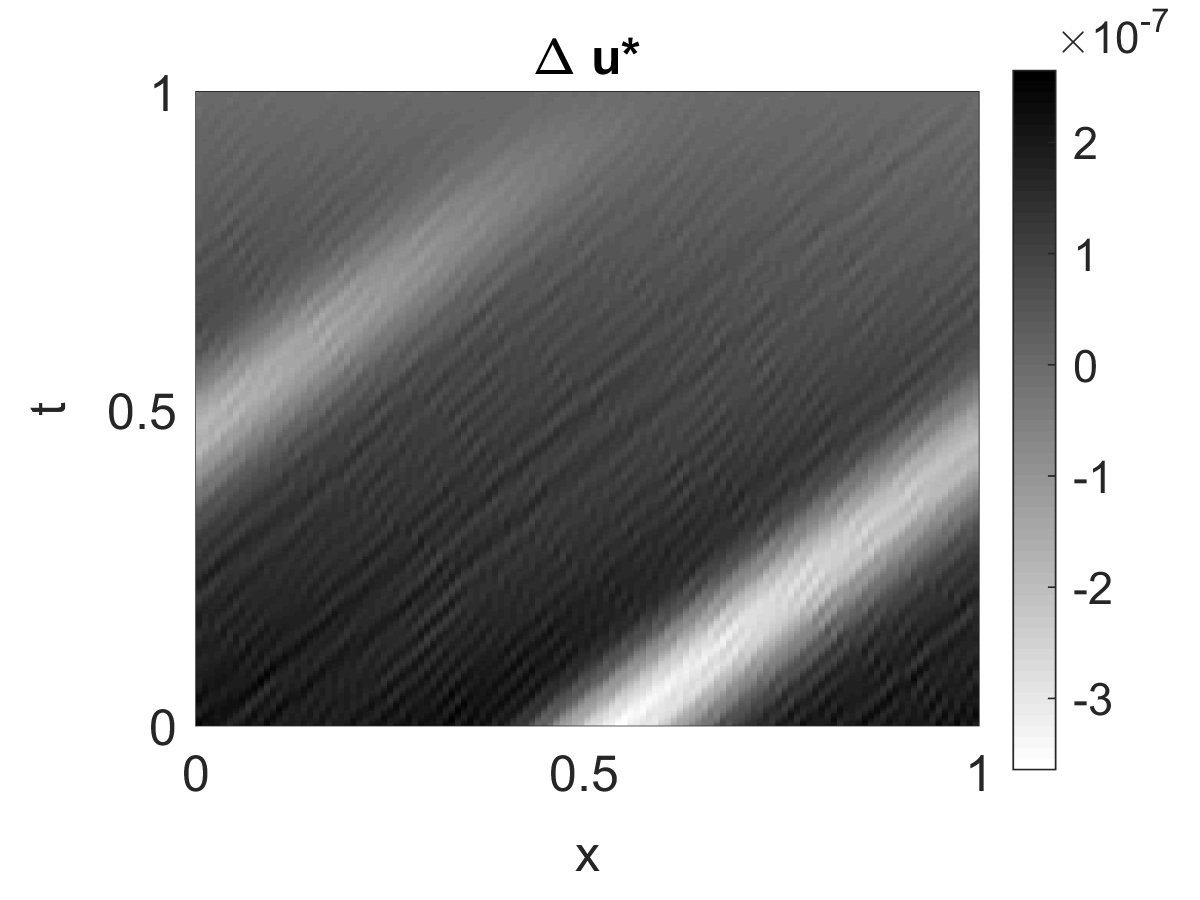}
	\caption{B1 - constant primal flow  solution with $u(x,t) = 1/2$ (top-left) and the corresponding analytical adjoint solution $u^*_{\mathrm{analytical}}$ (top-right). Mode-based adjoint solution $u^*_{\mathrm{mode-based}}$ (bottom-left) and difference $\Delta u^* = u^*_{\mathrm{mode-based}} - u^*_{\mathrm{analytical}}$ (bottom-right).
	\label{fig_results_c1}
	}
\end{figure}

\paragraph{B1 - Constant base flow}
The first setup of the Burgers equation is given by the initial condition of the primal problem $u(x,t_0) = 1/2$, which is thereby  the solution for all time.   
Thus, the primal and the adjoint equation reduces to a simple, constant transport, as can be seen from \eqref{eq_burgers} and \eqref{eq_burgersAdjoint}.
Likewise the calculation plan is just the current state and the obtained result since by this both are represented by $V$. 
The error is of the order of the linearization step used for the reference solution and within the dynamic Arnoldi method, 
so that the difference is likely created by the numerical linearization, see Fig.~\ref{fig_results_c1}.

\begin{figure}
  \centering
  \includegraphics[width = .49\textwidth]{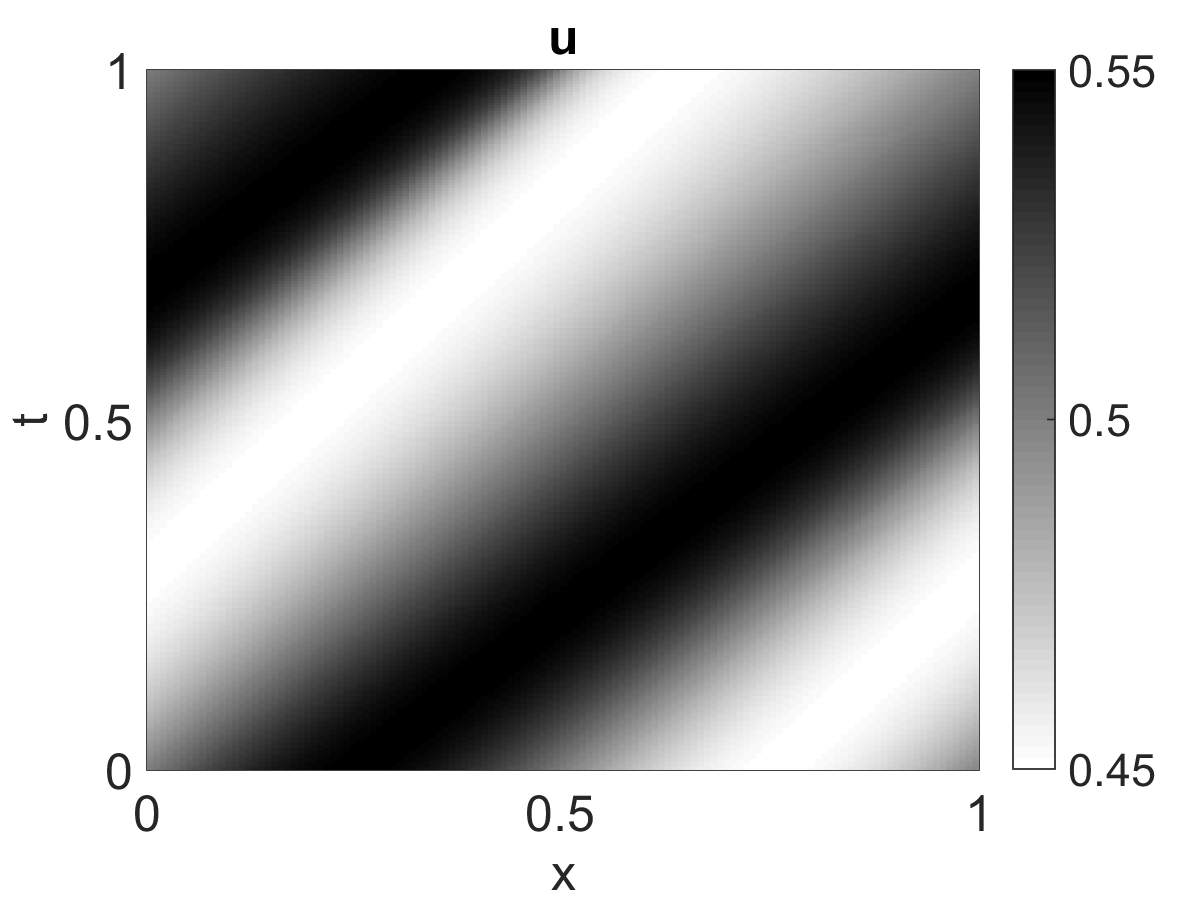}
  \includegraphics[width = .49\textwidth]{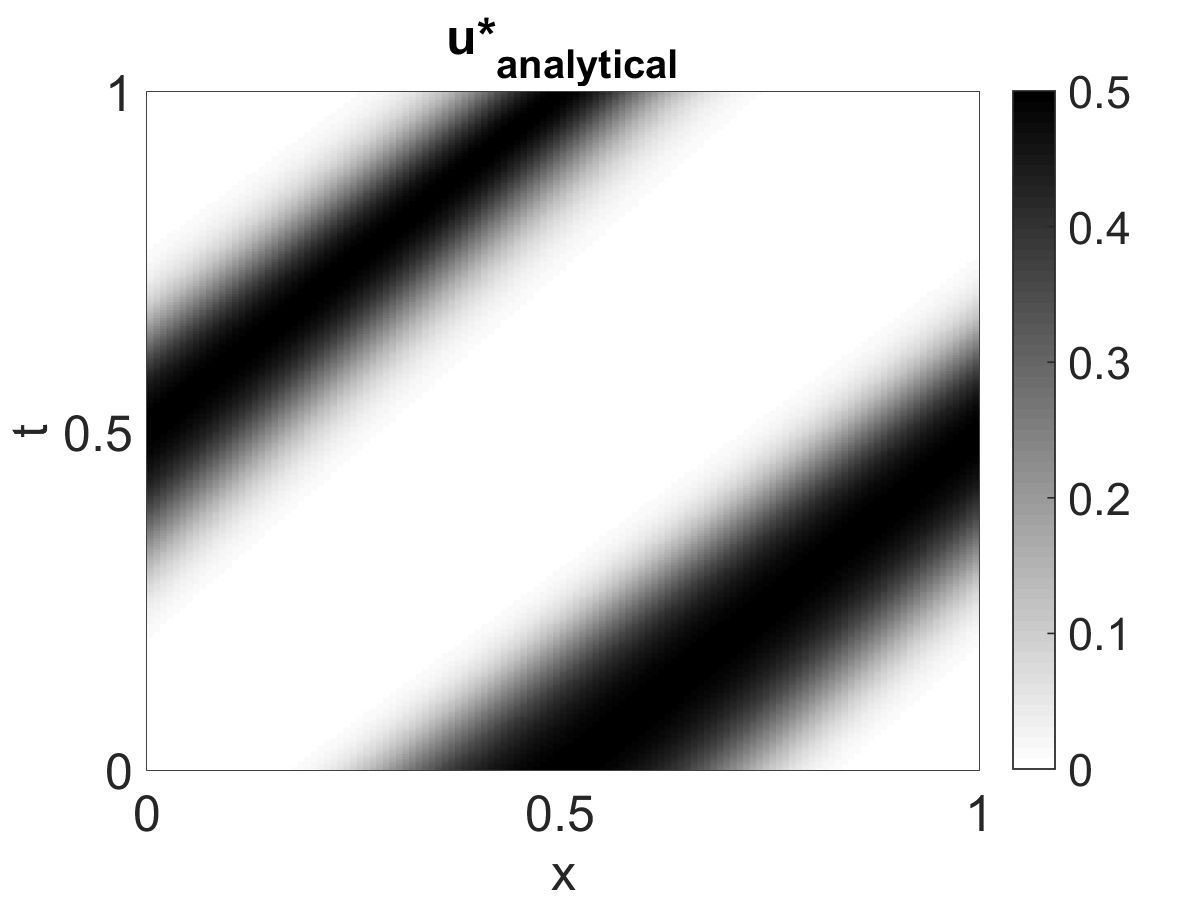}\\
  \includegraphics[width = .49\textwidth]{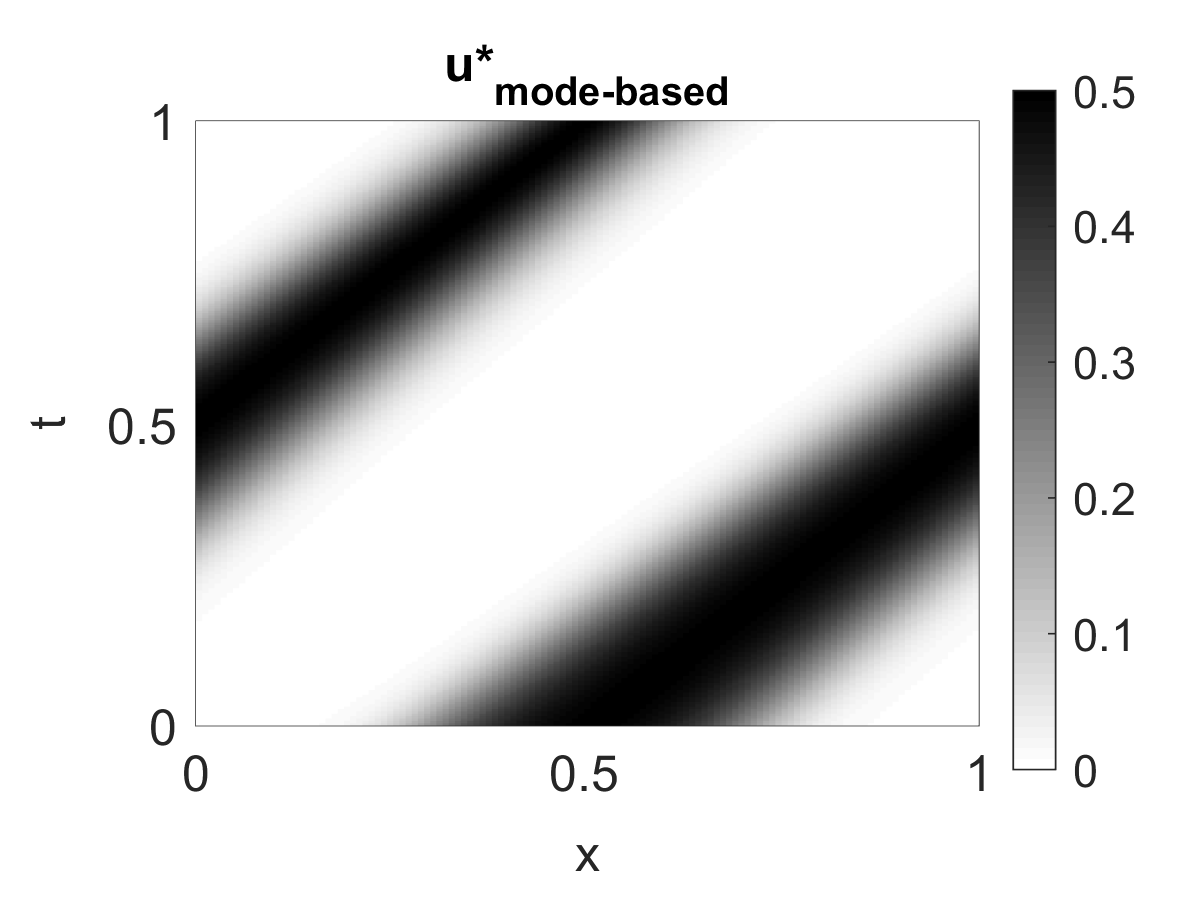}
  \includegraphics[width = .49\textwidth]{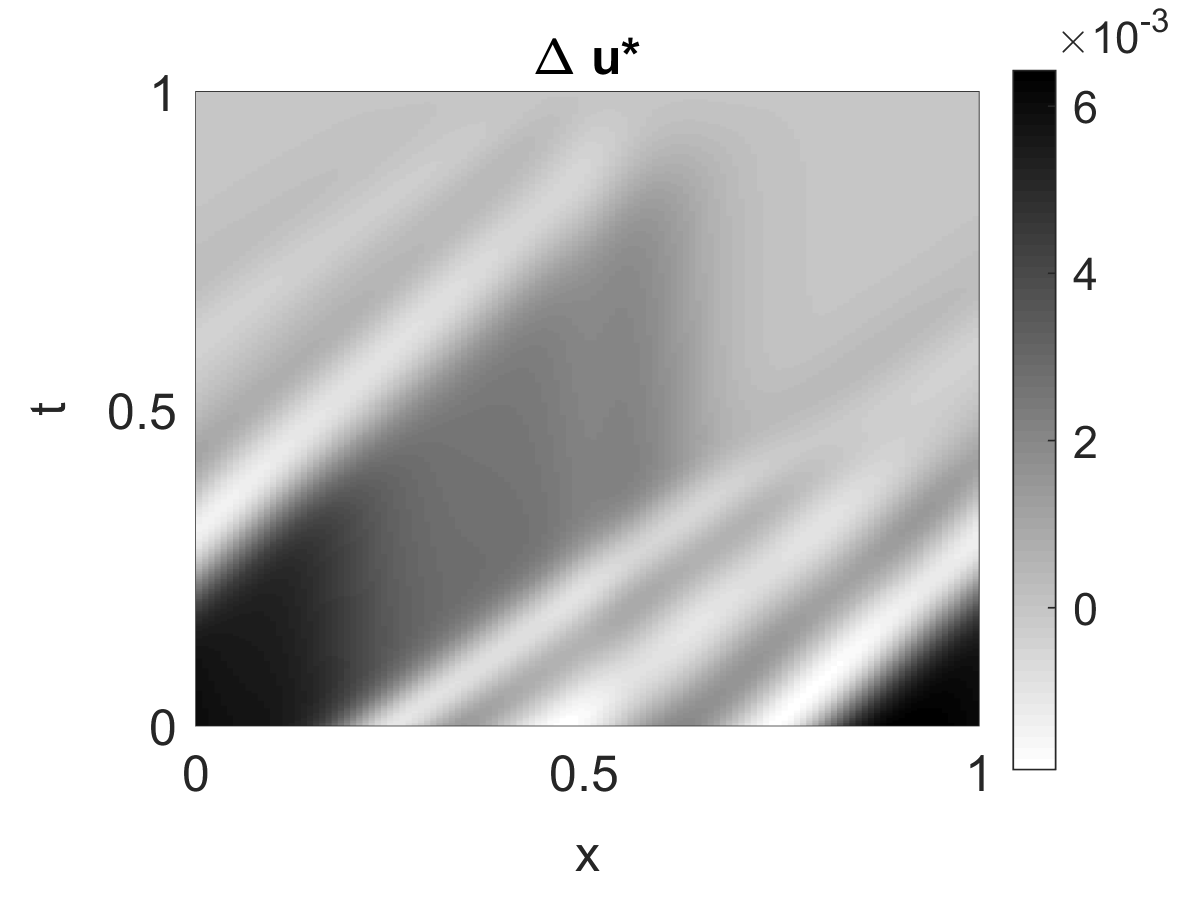}
  \caption{B2 - unsteady primal flow solution (top-left) and the corresponding analytical adjoint solution $u^*_{\mathrm{analytical}}$ (top-right). Mode-based adjoint solution $u^*_{\mathrm{mode-based}}$ (bottom-left) and difference $\Delta u^* = u^*_{\mathrm{mode-based}} - u^*_{\mathrm{analytical}}$ (bottom-right).    
  \label{fig_results_c2}
  }
\end{figure}

\paragraph{B2 - Unsteady base flow}
For this setup the primal initial condition is chosen to $u(x,t_0) = 1/2 + 1/20 \sin(x)$ leading to an unsteady flow solution.
The calculation plan is slightly longer and includes the right-hand-side (RHS) of the adjoint of the previous step $t_{i+1}$, calculated before.
This is little surprising since this changes only little from time step to time step and helps to construct the new RHS for time step $t_i$.
However, a larger deviation between mode-based and analytical solution is found.
The relative error is about 1\% with respect to the analytical reference solution, see Fig.~\ref{fig_results_c2}

Note, that for the first time step of the adjoint computation, where no previously computed RHS is available another calculation plan is needed.
However, this is a minor problem since it is possible to find a larger plan, which is sufficient.
The increased demand on computational time is negligible as this plan is just needed once for the first adjoint time step.

\begin{figure}
  \centering
  \includegraphics[width = .49\textwidth]{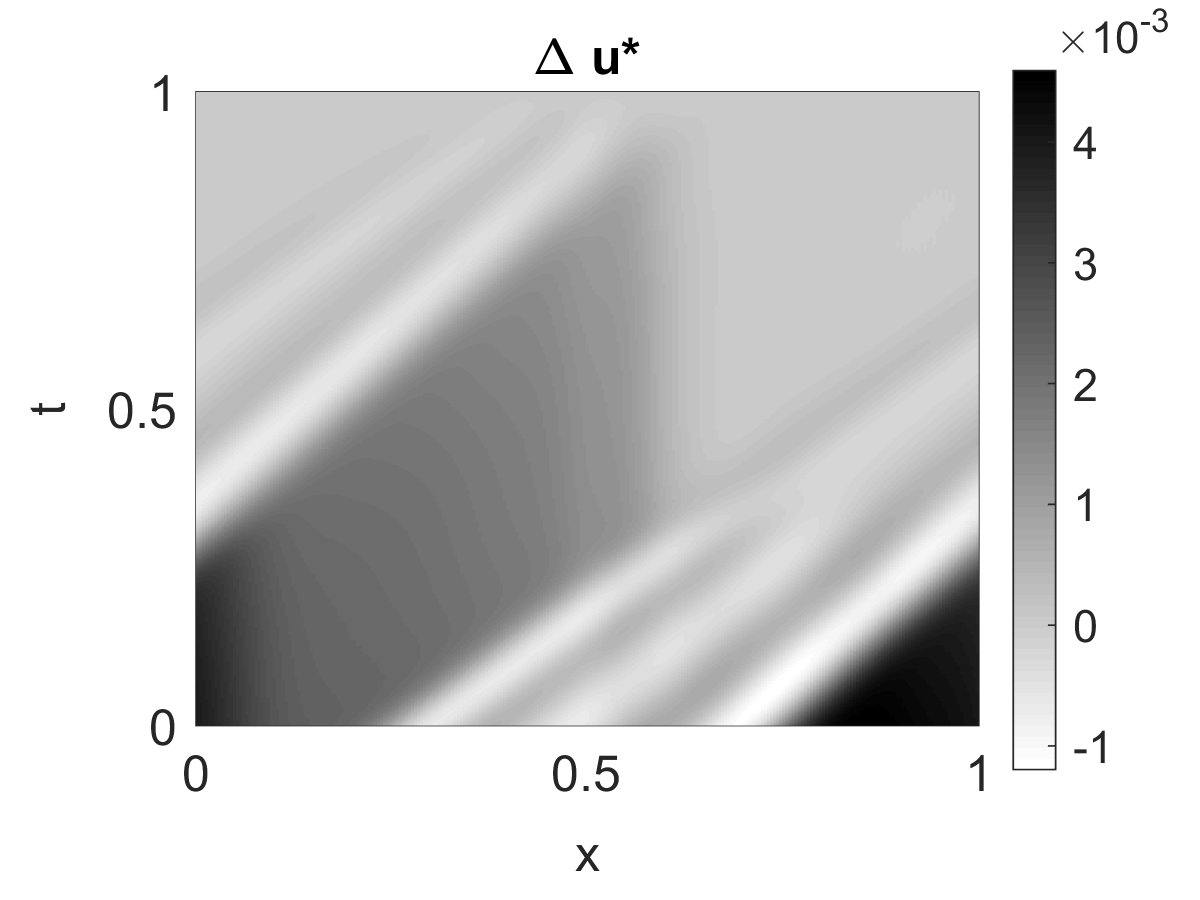}
  \includegraphics[width = .49\textwidth]{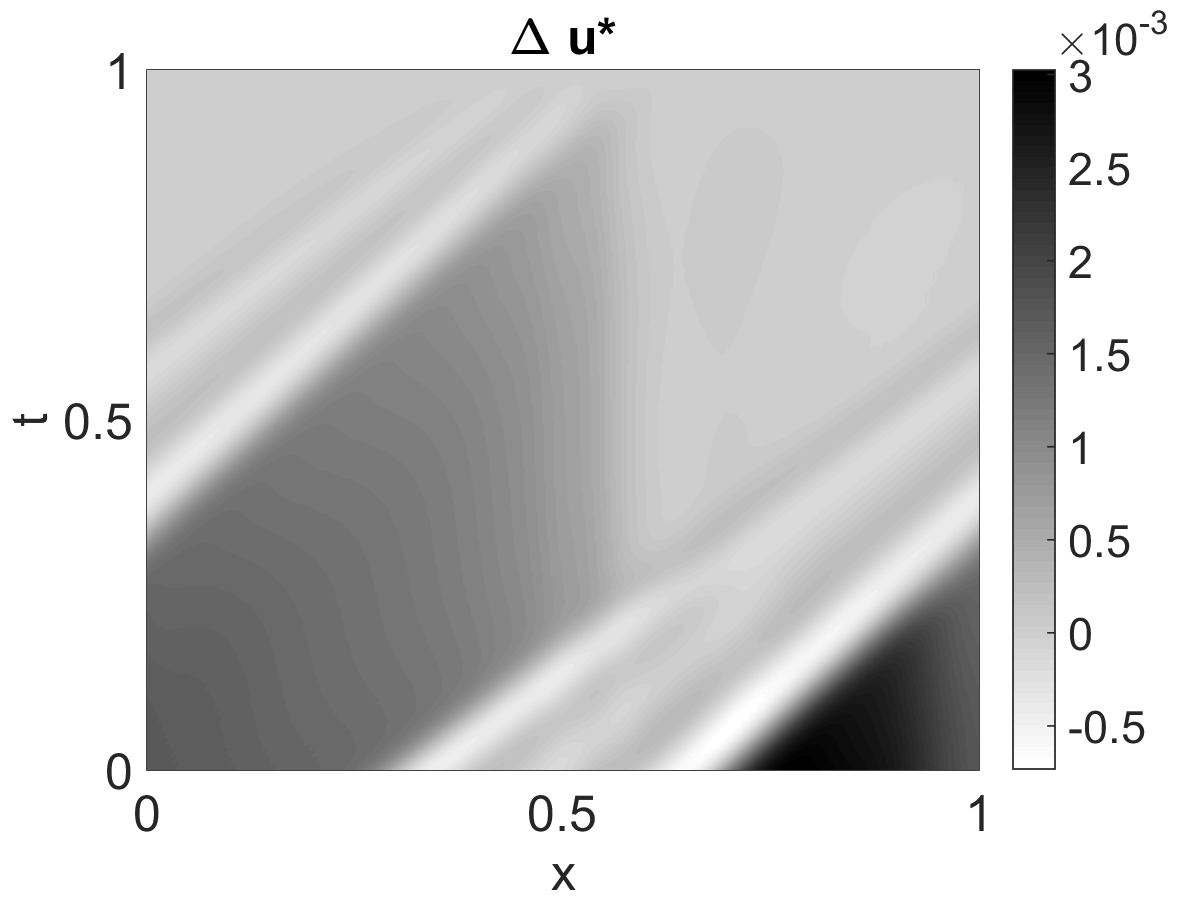}
  \caption{B2 - 
  Difference between mode-based and analytical adjoint solution $\Delta u^* = u^*_{\mathrm{mode-based}} - u^*_{\mathrm{analytical}}$ for setup B2 using $2 \times 128$ (left) and $4 \times 128$ (right) grid points for the spatial discretization.
  \label{fig_results_c2_2x_4x}
  }
\end{figure}

It is also important to note that the procedure is independent of the used spatial resolution.
Using the same calculation plan as before basically the same errors are found for double and quadruple the spatial resolution, see Fig.~\ref{fig_results_c2_2x_4x}.

\begin{figure}
  \centering
  \includegraphics[width = .49\textwidth]{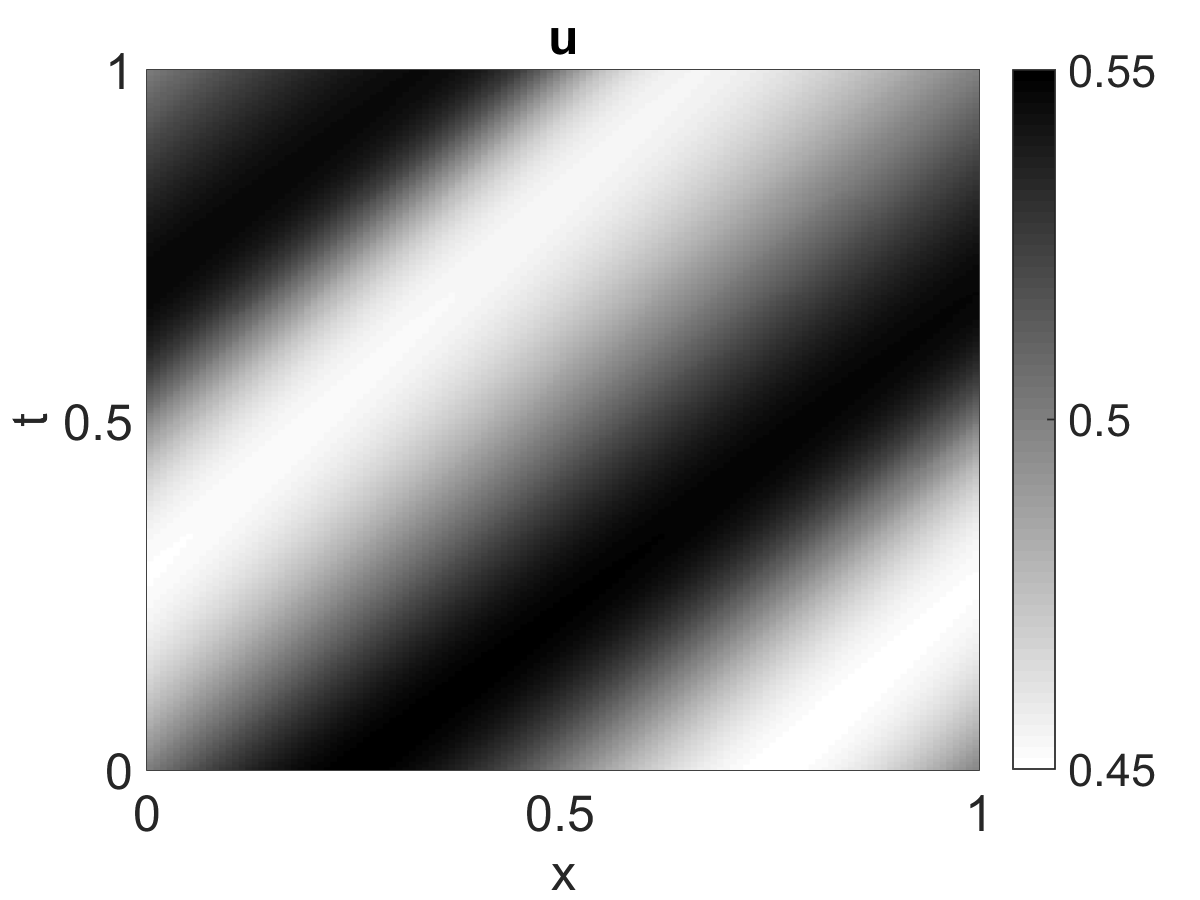}
  \includegraphics[width = .49\textwidth]{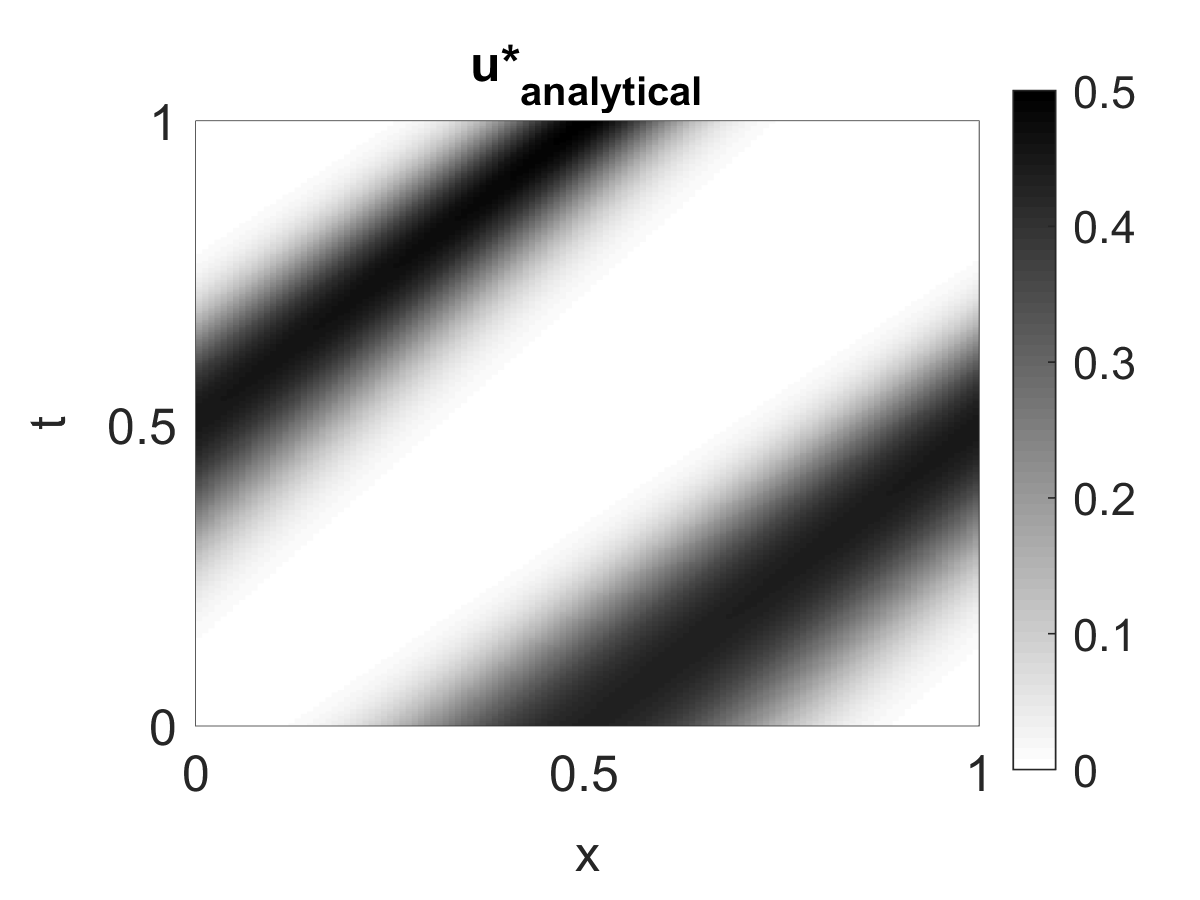}\\
  \includegraphics[width = .49\textwidth]{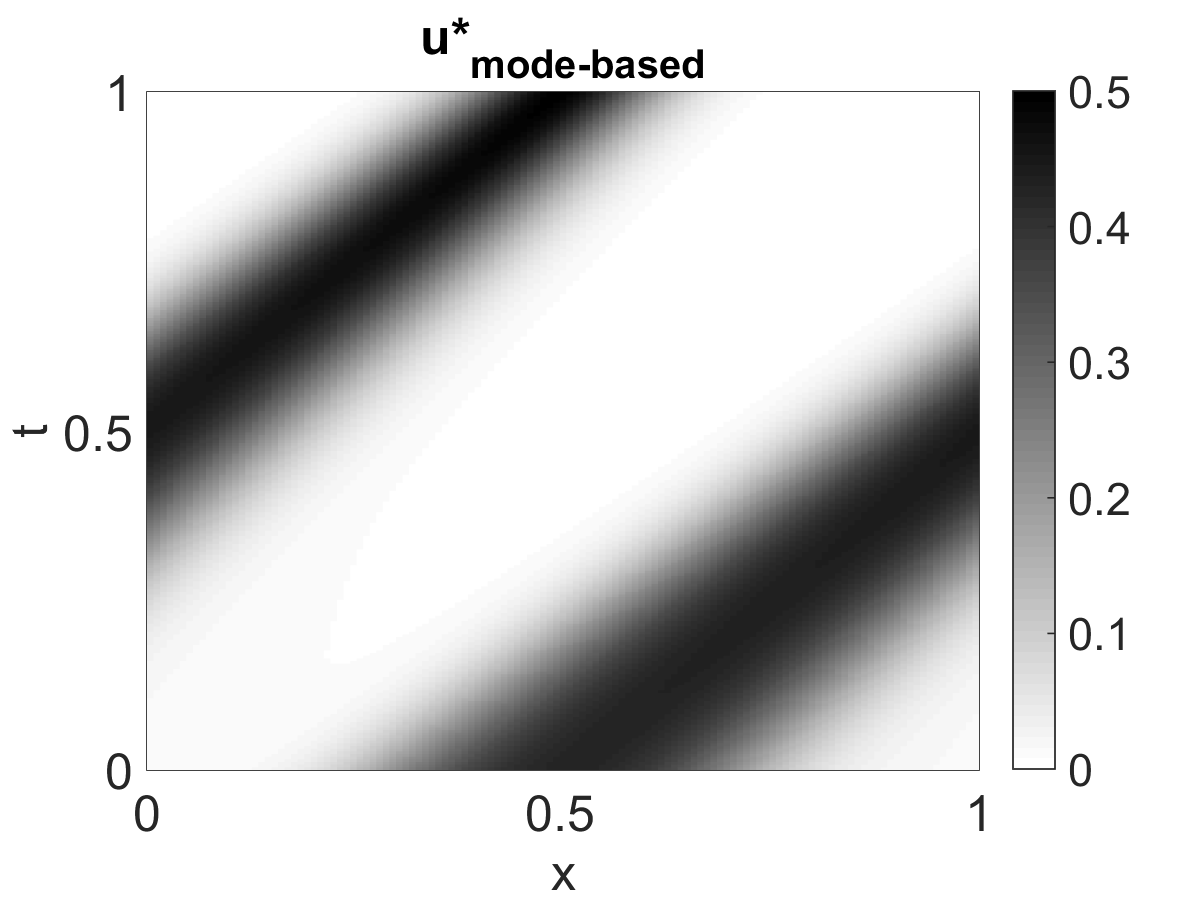}
  \includegraphics[width = .49\textwidth]{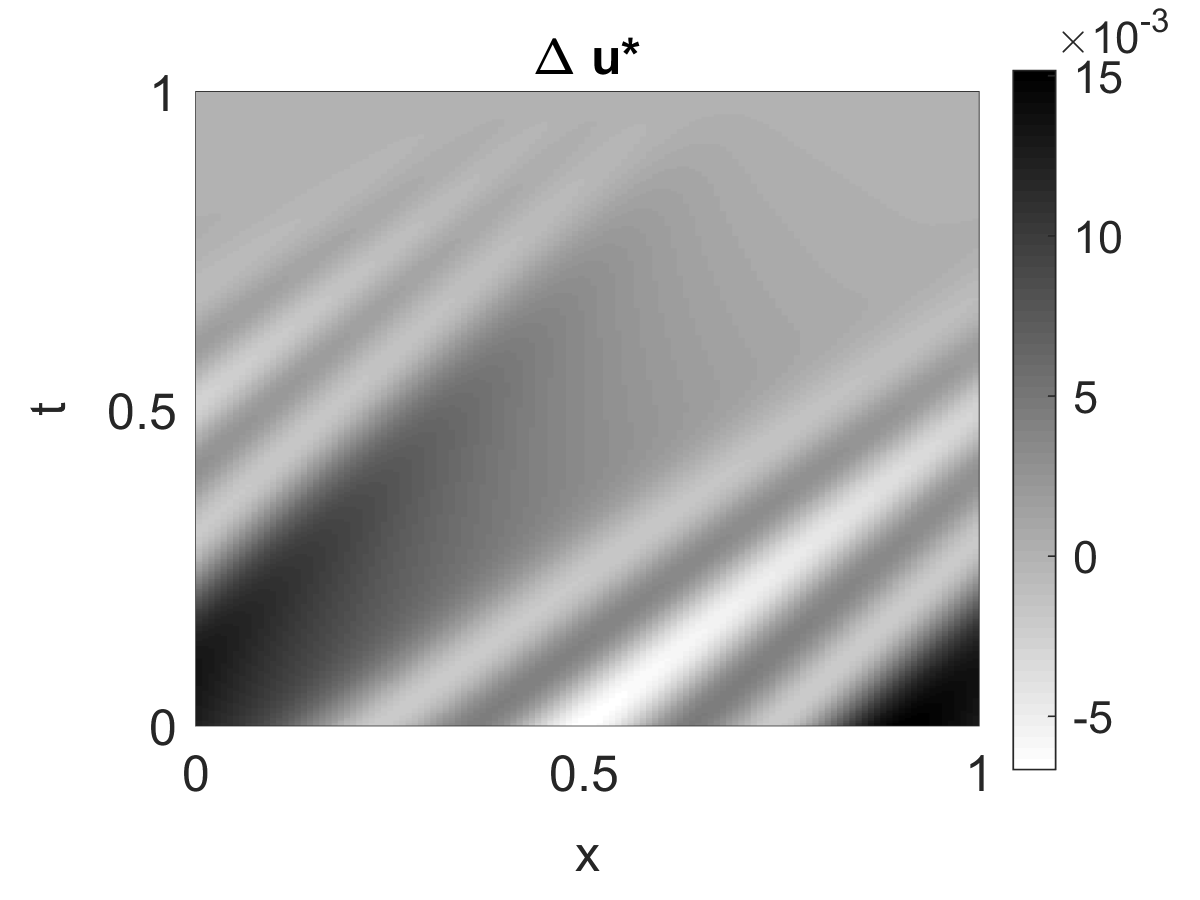}
  \caption{B3 - unsteady primal flow solution with friction (top-left) and the corresponding analytical adjoint solution $u^*_{\mathrm{analytical}}$ (top-right). Mode-based adjoint solution $u^*_{\mathrm{mode-based}}$ (bottom-left) and difference $\Delta u^* = u^*_{\mathrm{mode-based}} - u^*_{\mathrm{analytical}}$ (bottom-right).
  \label{fig_results_c3}
  }
\end{figure}

\paragraph{B3 - Unsteady base flow with friction}
This setup is dedicated to show the applicability of the mode-based adjoint method if the considered equation includes a second derivative $\partial_x^2$, e.g.~the friction term in \eqref{eq_burgers}.
Here, the friction constant is chosen to $\mu = 7.5\cdot 10^{-3}$.
The principal problem is that the corresponding friction part of the operator keeps his sign in the analytical adjoint equation, while the transport term changes.
Thus, the dynamic Arnoldi needs to separate the actions of the transport and the friction part.  
If this is not possible with an acceptable number of modes the approximation of adjoint operator is poor.
A suitable alternative to handle such a case is to split the equation and treat the transport terms and the friction term individually, by the same calculation plan.
Both resulting adjoint operators are simply added in order to construct the complete adjoint operator.
This procedure leads to errors on the level of case (B2), see Fig.~\ref{fig_results_c3}.
%


\subsection{\textsc{Euler} Equations \label{ResultEuler}}
In the following the more complex problem of coupled equations is discussed on the base of the one dimensional Euler equations
\begin{eqnarray}
  \partial_t \rho  	+ \partial_x (\rho u )    					&=& 0 \no\\ 
  \partial_t \rho u 	+ \partial_x (\rho u u ) 	+ \partial_x p  		&=& 0 \no\\ 
  \partial_t  p  	+ \gamma \partial_x (p  u ) 	- (\gamma -1) u \partial_x p 	&=& 0. \label{eq_euler}
\end{eqnarray}
Therein $\rho$ denotes the density, $u$ the velocity, $p$ the pressure and $\gamma$ the adiabatic exponent which is assumed to be $1.4$.
The equations are discretized again by means of a finite difference approach in space.
The computational domain of length $2\pi$ is resolved by $128$ equidistantly distributed and periodic if not stated otherwise.
Again fourth order differentiation schemes are employed.
The computational time span $t_0$ to $t_{\mathrm{end}}$ is separated into $171$ time steps using of CFL-condition of $0.75$, based on the base flow velocity plus the speed of sound.
A fourth order Runge-Kutta scheme is employed for the time-wise integration.
Again, the Arnoldi-based adjoint and the reference adjoint are using the same discretization.

The corresponding adjoint equations of \eqref{eq_euler} are derived and discussed in \cite{Lemke2015}.
All adjoint computations, based on the proposed method or on an analytical derivation, are initialized by means of a Gaussian disturbance in $p^*$
\begin{equation}
 p^*(x,t_{\mathrm{end}}) = 5 \cdot \exp\left(-\frac{(x - x_0)^2}{(10 \Delta x)^2} \right),
\end{equation}
with $x_0$ as center of the computational domain and the grid spacing $\Delta x$, at the end of the computational time $t_{\mathrm{end}}$.

For all setups the employed calculation plans, required for the mode-based adjoint method, can be found in Tbl.~\ref{app_tbl_calcplans_euler} in the appendix.
Within this calculation plans the RHS of the previous time-step was not incorporated, since this produced smaller training plans, see discussion in Sec.~\ref{sec_training_of_the_method}. 

To allow a discussion of the calculation plans a linearized version of \eqref{eq_euler} is derived as
\begin{equation}
\partial_t \begin{pmatrix}  \delta \rho \\ \delta u \\\delta p \end{pmatrix}  
+
\partial_x \left(
\underbrace{ 
\begin{pmatrix} 
  u_0 & \rho_0 & 0        \\ 
  0   & u_0    & 1/\rho_0 \\
  0   & \gamma p_0 & u_0 
\end{pmatrix}}_{=\mathcal{A}} 
\begin{pmatrix} \delta \rho \\\delta u \\\delta p \end{pmatrix}  
\right)
= 
0 . \label{eq_lin_euler}
\end{equation}
under the assumption of spatial  constant base flows ($\partial_x \rho_0 = \partial_x u_0 = \partial_ xp_0 = 0  $).
Note that $\mathcal{A}$ is not the desired operator $A$, which includes spatial discretization and possible boundary treatment.

\paragraph{E1 - No base flow} 
For the first setup a flow at rest condition with $\rho(x,t_0) = 1$, $u(x,t_0) = 0$ and $p(x,t) = 1.5$ is chosen.
Thus, the Euler equations reduce to purely acoustic equations.
One might expect that this is trivial, since the acoustic equations are known to be (with an additional rescaling) self-adjoint. 
However, the structure of $\mathcal{A}$
\begin{equation}
  \mathcal{A}  = 
  \begin{pmatrix} 
    0 & \rho_0 & 0        \\ 
    0   & 0    & 1/\rho_0 \\
    0   & \gamma p_0 & 0 
  \end{pmatrix} 
  \qquad 
  \mathcal{A}^T  = 
  \begin{pmatrix} 
    0 & 0 & 0        \\ 
    \rho_0  & 0    & \gamma p_0 \\
    0   & 1/\rho_0 & 0 
  \end{pmatrix}
  \label{eq_lin_euler_no_base}
\end{equation}
reveals, that $\delta \rho $ is driven by $\delta u$, but the opposite does not hold.
The adjoint equation of \eqref{eq_lin_euler_no_base} is given by 
\begin{equation}
  \partial_t q^* + \partial_x \mathcal{A}^T q^* = \partial_t q^* + \mathcal{A}^T \partial_x q^* = 0      
\end{equation} 
with $q^* = \left[\rho^*, u^*, p^* \right]^T $ as adjoint variable and using that $q_0 = [\rho_0, u_0, p_0]$ is constant in space. 
The matrix $\mathcal{A}^T$ has zero first line, reflecting that changes in $\rho$ have no effect on the other quantities whatsoever.
This reveals an important structural feature of systems with multiple coupled equations: one variable influences a second one while there is no (or a structurally different) feedback.
This poses a severe difficulty for the classical Arnoldi or block-Arnoldi.
As can be observed easily, the vectors created by iterating the matrix $A$ produce, the same entries for $\rho$ and $u$ up to a factor.
Applying the Arnoldi directly to $\mathcal{A}^T$ would create structurally different vectors with a zero entry for $\rho$.
Therefore, it is in general not possible to span the Krylov space of the adjoint operator by Krylov vectors\footnote{Neither the Krylov space of $A$ nor the one of $A^T$ span the full discrete vector space, at least for exact arithmetic.} of $\mathcal{A}$.

To circumvent this problem we allow to modify the vectors in the dynamic Arnoldi, which was already discussed in Sec.~\ref{sec_dam}.
For this case, the calculation plan employs the actual adjoint state vector and three modified vectors, which results from application of the primal right-hand-side, see Tbl.~\ref{app_tbl_calcplans_euler}. 
It can by easily checked, that for this simple example, the calculation plan is able to reproduce the analytical solution.

\begin{figure}
  \centering
  \includegraphics[width = .49\textwidth]{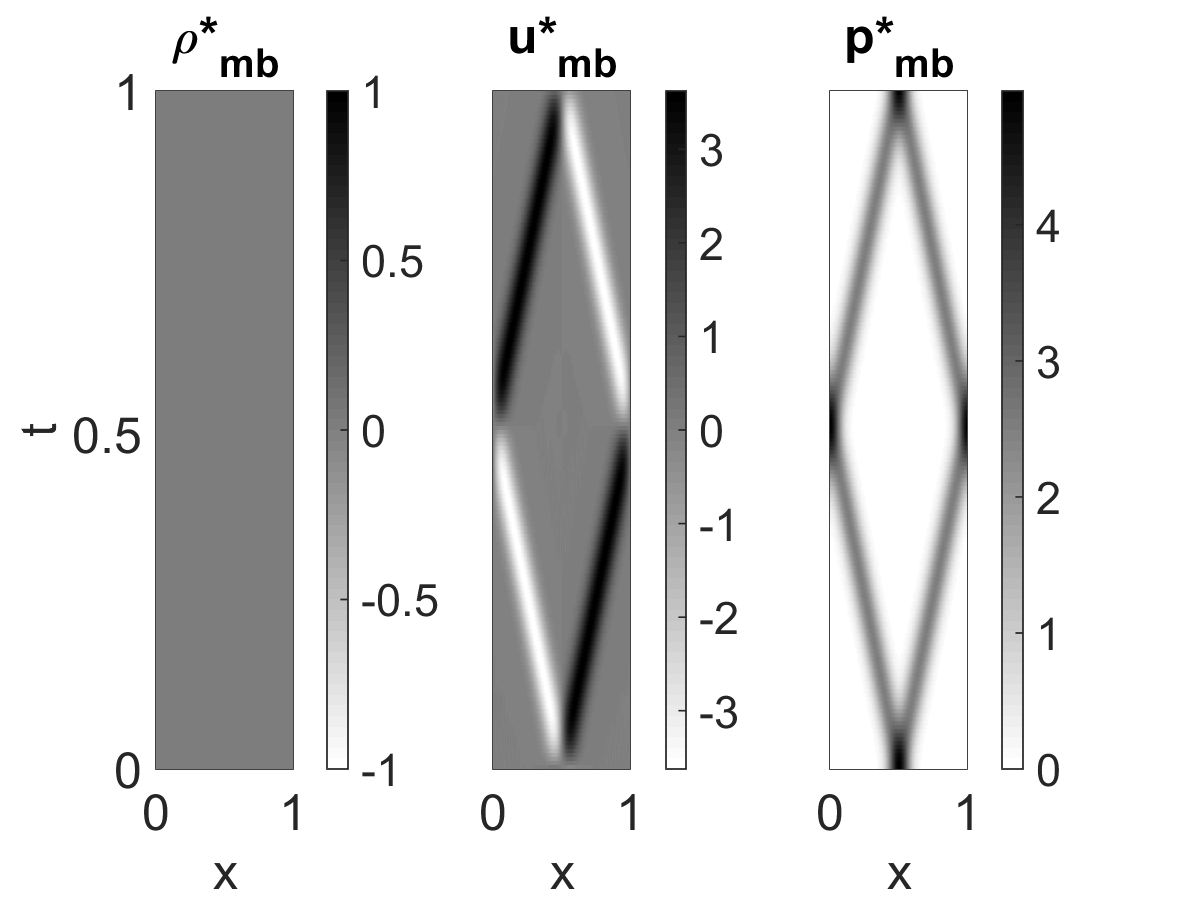}
  \includegraphics[width = .49\textwidth]{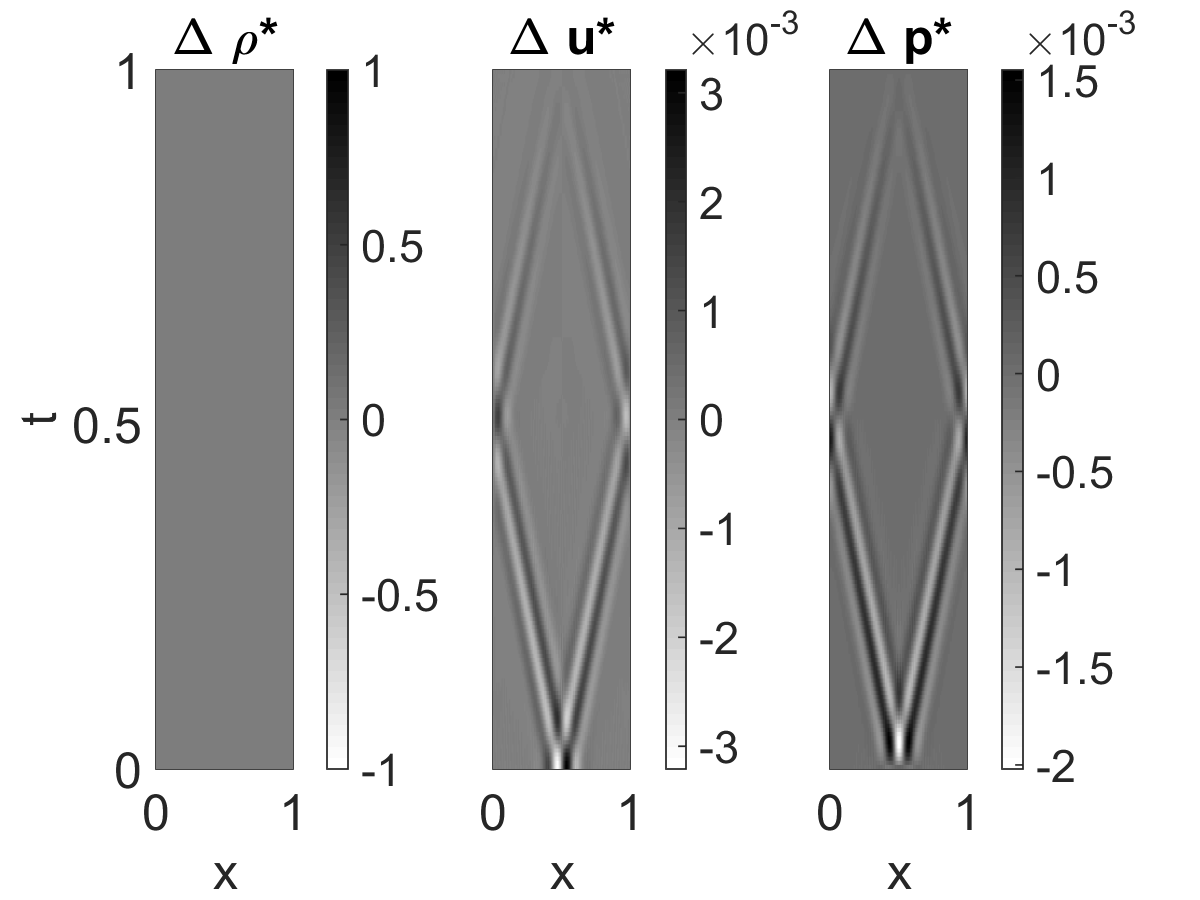}
  \caption{E1 - No base flow condition. Mode-based (mb) adjoint solution $u^*_{\mathrm{mb}}$ (left) and difference $\Delta q^* = q^*_{\mathrm{mb}} - q^*_{\mathrm{analytical}}$ (right) for all quantities.  
  \label{fig_results_c4}
  }
\end{figure}

The numerical results of the adjoint and the deviation with respect to the analytical solution are shown in Fig.~\ref{fig_results_c4}.
The adjoint solution is defined by acoustic characteristics in $u^*$ and $p^*$, while $\rho^*$ remains zero, as expected from the analytical solution.
The error is  similar to the cases related to the Burgers equation.

\paragraph{E2 - No base flow - high pressure} 
The primal computational variables are usually of very different magnitude \cite{SesterhennMuellerThomann1999}, which can be problematic for the error control of iterative procedures and linearization.  
This setup repeats the former case with physical units, leading to adjoint variables of different magnitudes.  
The initial pressure is $p(x,t_0) = 10^5$ while the density is $\rho(x,t_0) = 1$.

\begin{figure}
  \centering
  \includegraphics[width = .49\textwidth]{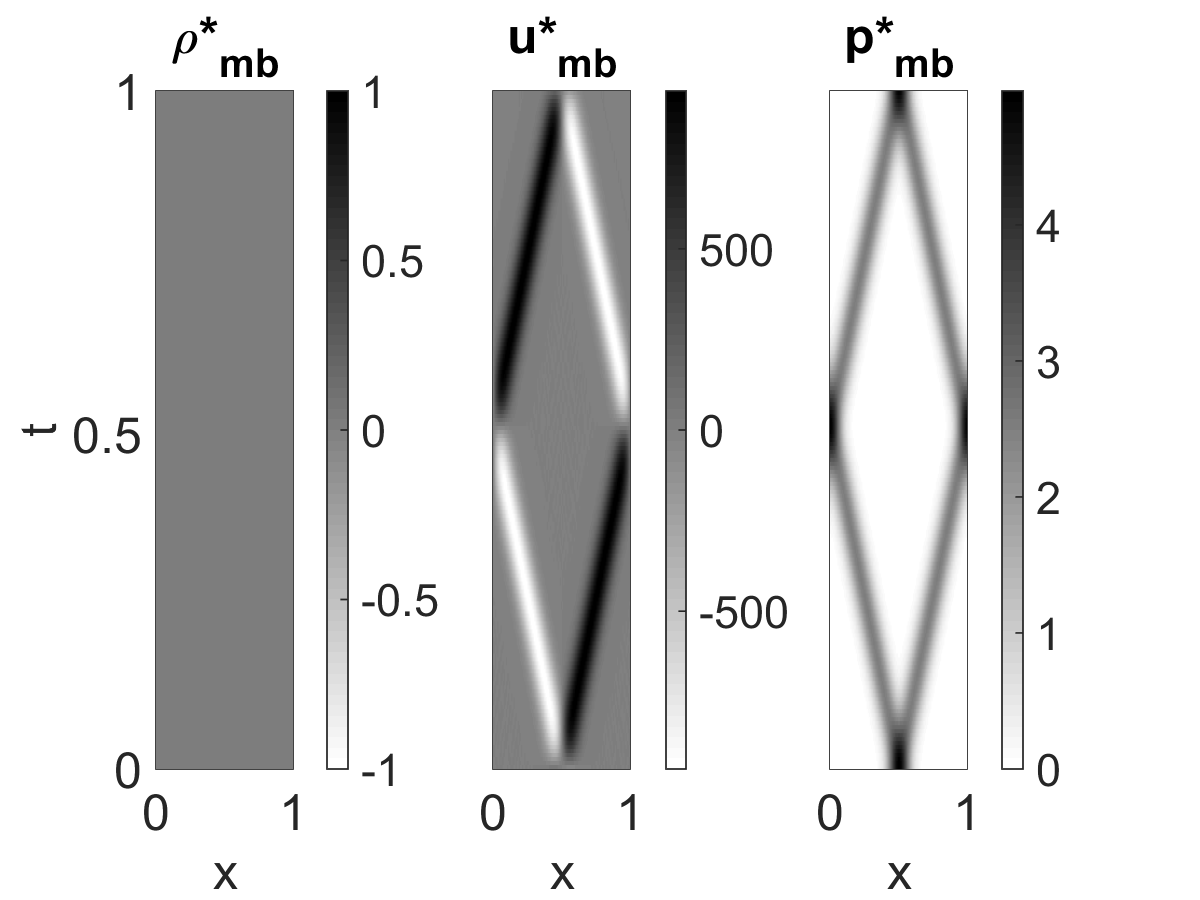}
  \includegraphics[width = .49\textwidth]{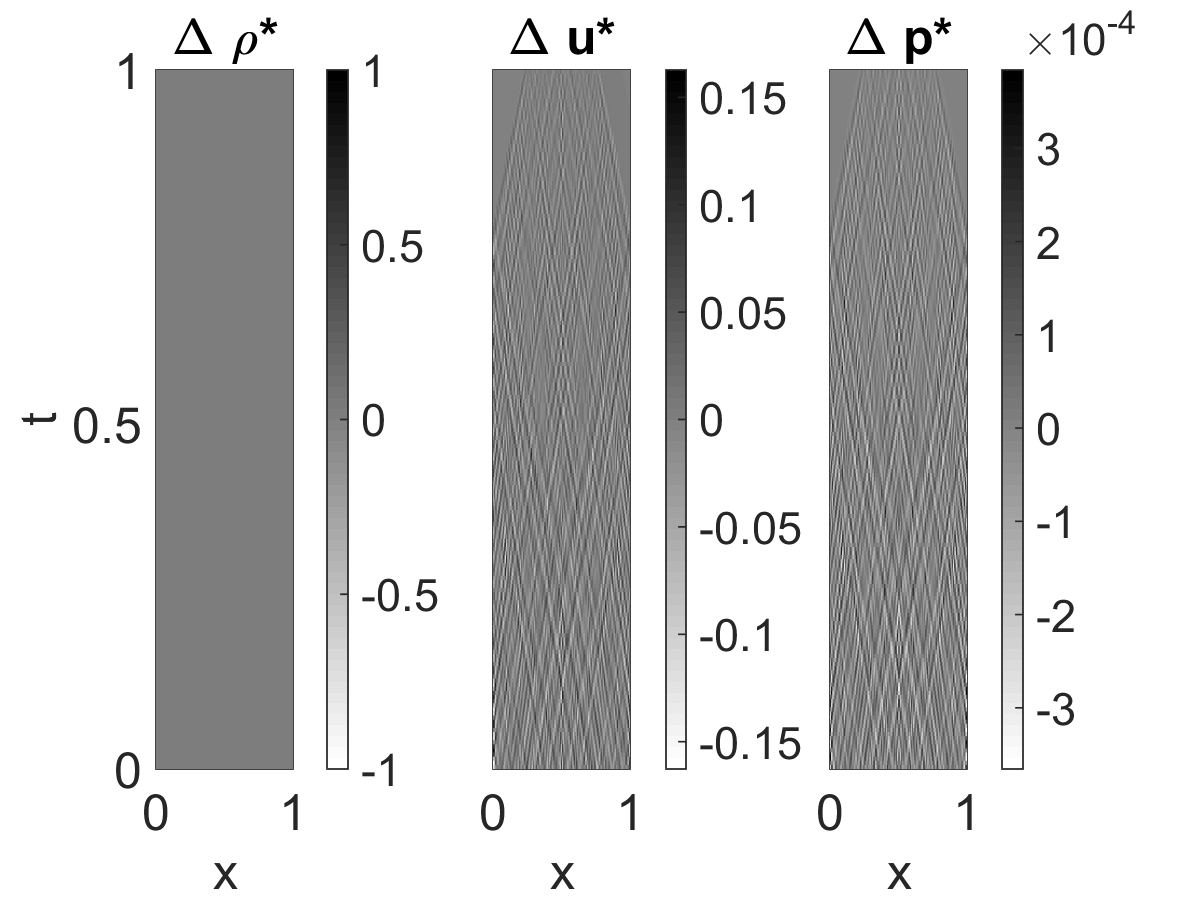}
  \caption{E2 - No base flow condition with $p_0=10^5$. Mode-based (mb) adjoint solution $u^*_{\mathrm{mb}}$ (left) and difference $\Delta q^* = q^*_{\mathrm{mb}} - q^*_{\mathrm{analytical}}$ (right) for all quantities. Please note the different magnitudes of the adjoint quantities.
  \label{fig_results_c5}
  }
\end{figure}

A slightly different calculation plan was found by the training, see Tbl.~\ref{app_tbl_calcplans_euler}.
The results match those of the previous case in terms of accuracy, see Fig.~\ref{fig_results_c5}.
If the calculation plan of (E1) is used the resulting adjoint solution is generally consistent with the analytical solution, but characterized by high-frequency fluctuations.

\begin{figure}
  \centering
  \includegraphics[width = .49\textwidth]{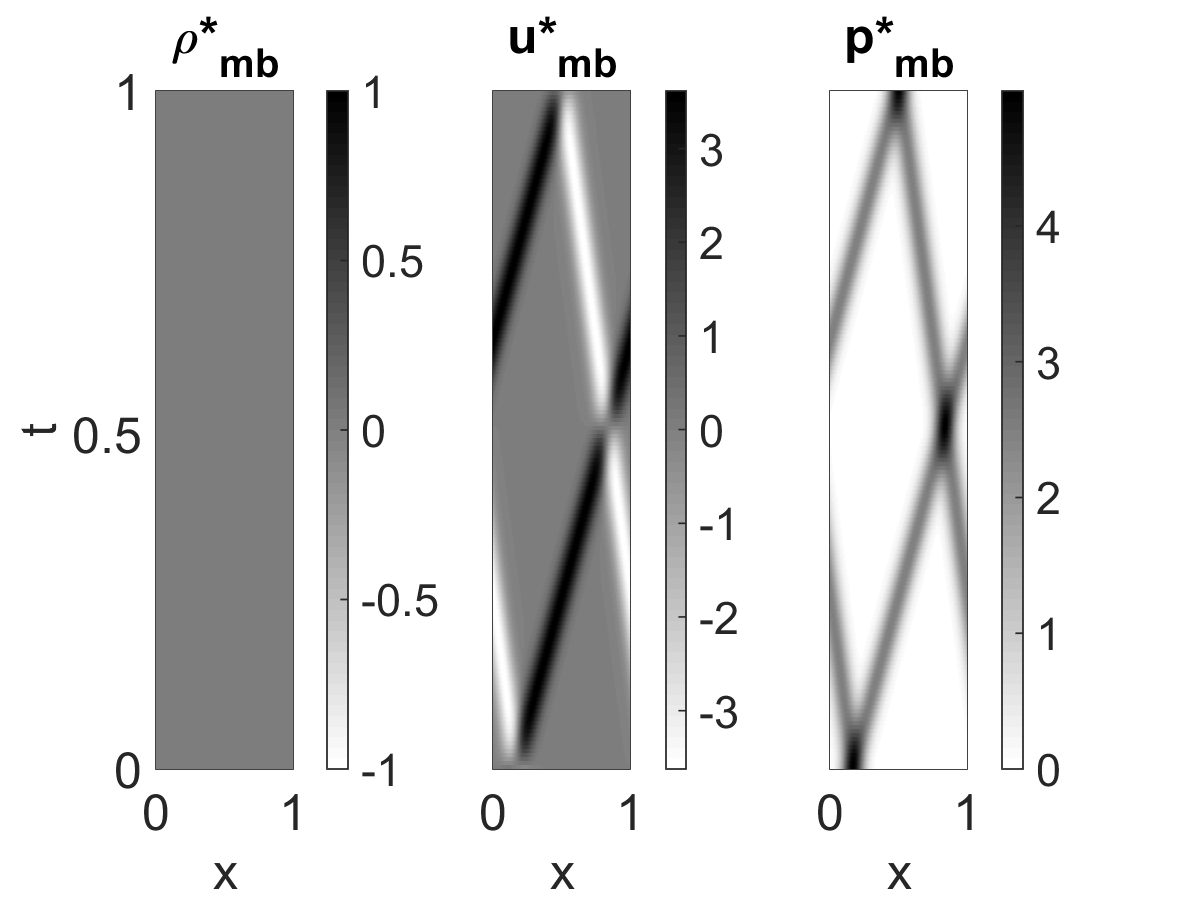}
  \includegraphics[width = .49\textwidth]{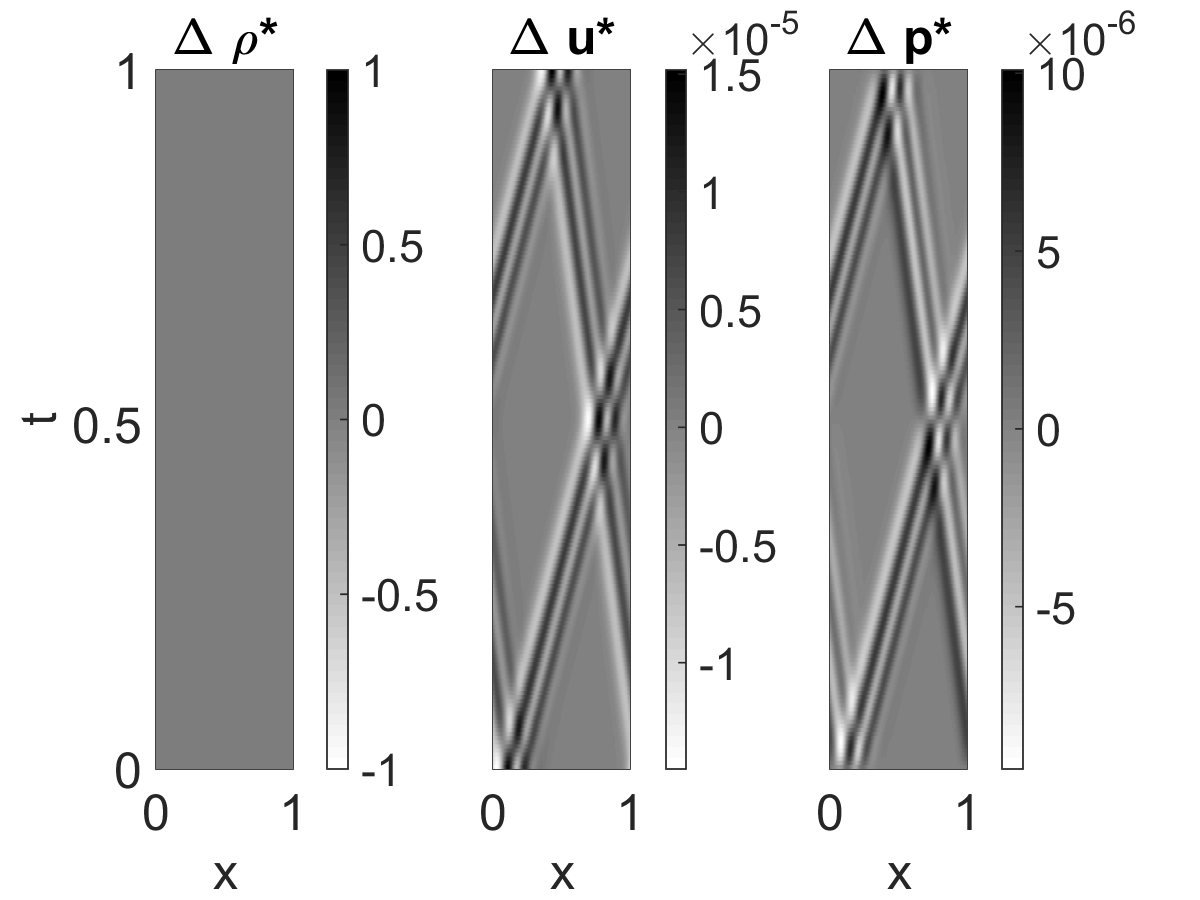}
  \caption{E3 - Steady base flow - intermediate Ma number. Mode-based (mb) adjoint solution $u^*_{\mathrm{mb}}$ (left) and difference $\Delta q^* = q^*_{\mathrm{mb}} - q^*_{\mathrm{analytical}}$ (right) for all quantities.
  Please note the skewness of the characteristics due to the presence of a base flow.
  \label{fig_results_c6}
  }
\end{figure}

\paragraph{E3 - Steady base flow - intermediate Ma number}
The primal initial conditions are defined by $\rho(x,t_0) = 1$, $p(x,t_0) = 1.5$ and $u(x,t_0) = 1/3c$, with $c = \sqrt{\gamma p_0/\rho_0}$ as the speed of sound.
Within this setup a steady non-zero base flow velocity, which breaks the self-adjoint structure of the governing system, is analyzed.
Despite of the base flow, the analytic solution yields $\rho^*(x,t) = 0$ for all time steps.
Considering the operator $\mathcal A$ in \eqref{eq_lin_euler} this is difficult to realize as $\delta \rho$, previously used solely for the representation of $u^*$ and $p^*$, acts on its own.
In order to remove this dependency and allow for a zero solution in $\rho^*$ a modified input vector based on the current adjoint state is employed.
Using only the pressure part of $q^*$ as additional vector for the dynamic Arnoldi the method is enabled to remove the entanglement as the wrong extra  part in the modes  can be be represented as difference, see Tbl.~\ref{app_tbl_calcplans_euler} for the calculation plan.
Using six base vectors a similar quality of the adjoint solution, which is characterized by a skewness of the characteristics due to the base flow, is found, see Fig.~\ref{fig_results_c6}.

\paragraph{E4 - No base flow - open boundaries}
For this setup the periodic boundary conditions are replaced by non-reflecting open boundaries, the spatial discretization is of second order.
In more detail, characteristic-boundary conditions \cite{Lemke2015,PoinsotLele1992} are combined with a quadratic sponge layer \cite{Mani2012} with acts on 10\% of the computational domain on both sides.
All other parameters of this setup corresponds to (E2).

\begin{figure}
  \centering
  \includegraphics[width = .49\textwidth]{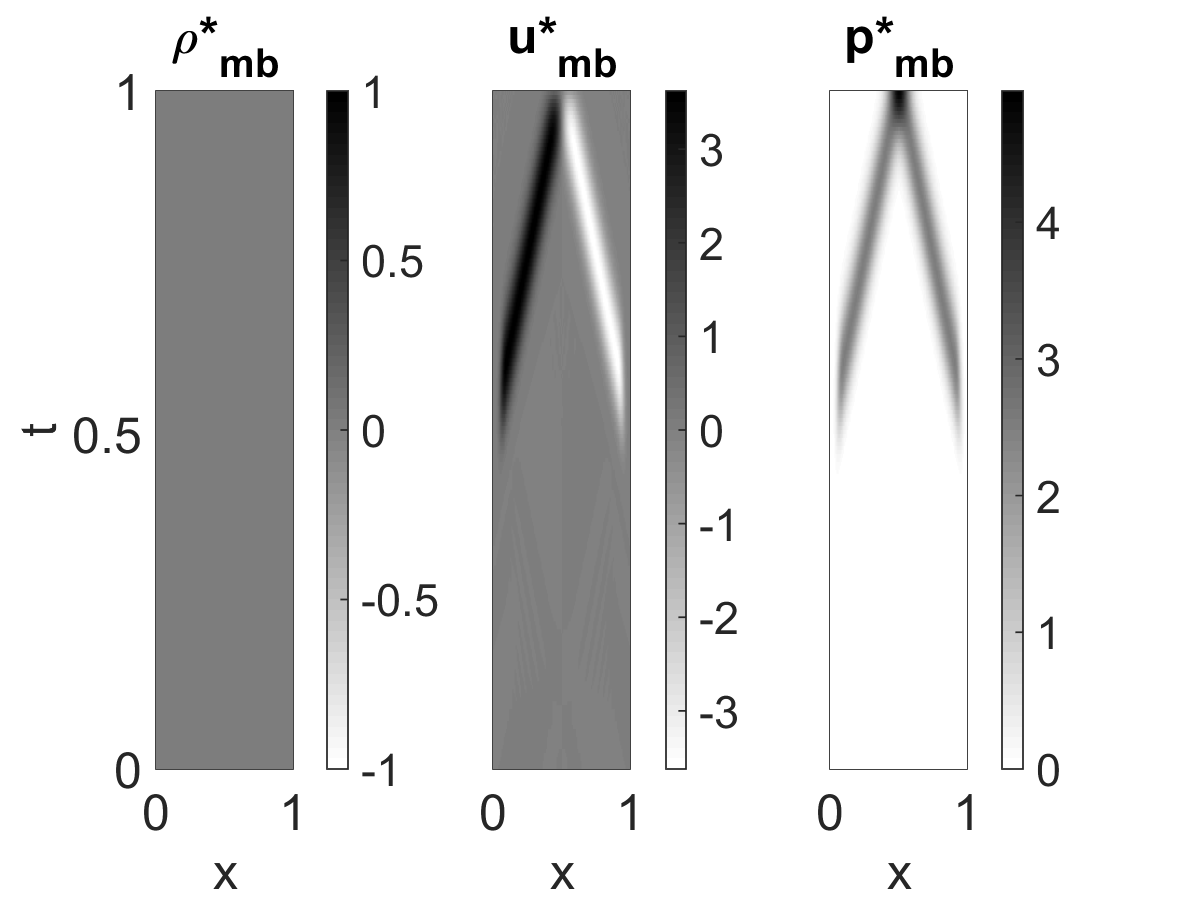}
  \includegraphics[width = .49\textwidth]{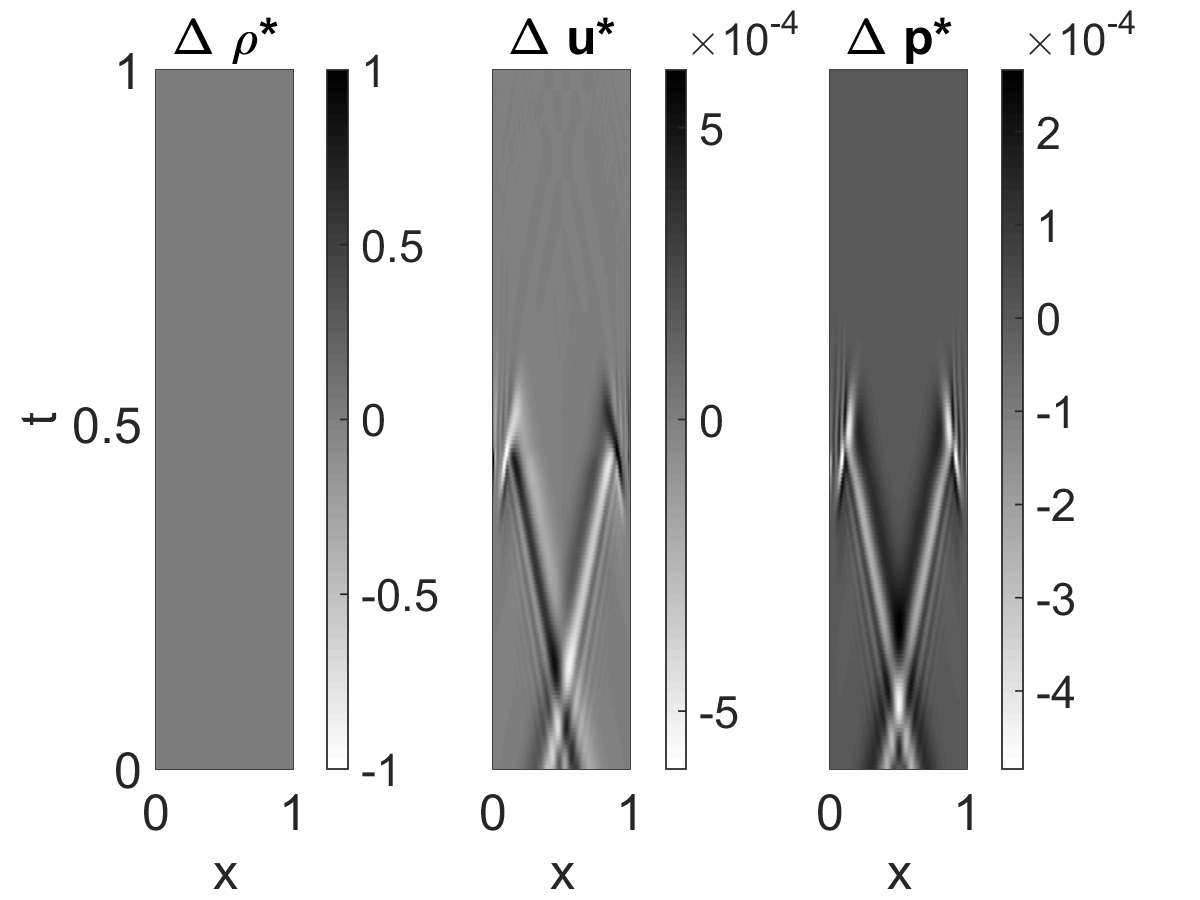}
  \caption{E4 - No base flow - open boundaries. Mode-based (mb) adjoint solution $u^*_{\mathrm{mb}}$ (left) and difference $\Delta q^* = q^*_{\mathrm{mb}} - q^*_{\mathrm{analytical}}$ (right) for all quantities.
  \label{fig_results_c7}
  }
\end{figure}

The presence of the damping sponge requires an additional number of modes.
The calculation plan results in $10$ calls of the primal right hand side, see Tbl.~\ref{app_tbl_calcplans_euler}.
Figure \ref{fig_results_c7} shows that the quality of the resulting adjoint solution is similar to case (E2), at least before the pulses reach the boundaries.
Thereafter, only slight reflections are found in contrast to the analytic solution, see \cite{Lemke2015} for details on the adjoint sponge layer.

Please note, that other boundary conditions, for example (no)-slip walls, might have to be treated individually.
As the present approach approximates the discrete adjoint operator of the governing equations (corresponding to the automatic difference method), the same problems are expected to arise.
A similar treatment should be possible.
Further discussion can be found in \cite{GilesDutaMuellerPierce2003,GilesPierce1997}.


\subsection{O1 - Noise Cancellation}
\begin{figure}
  \centering
  \includegraphics[width = .49\textwidth]{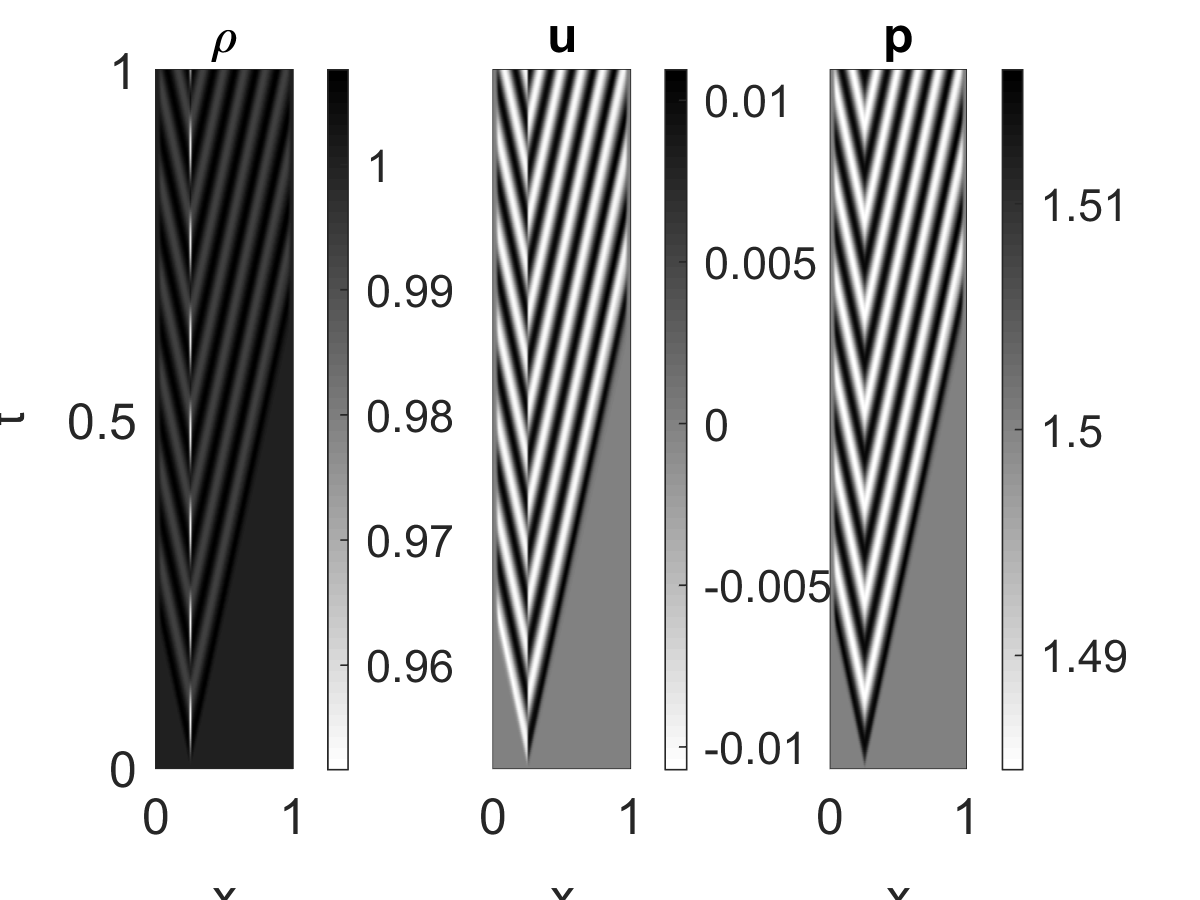}
  \includegraphics[width = .49\textwidth]{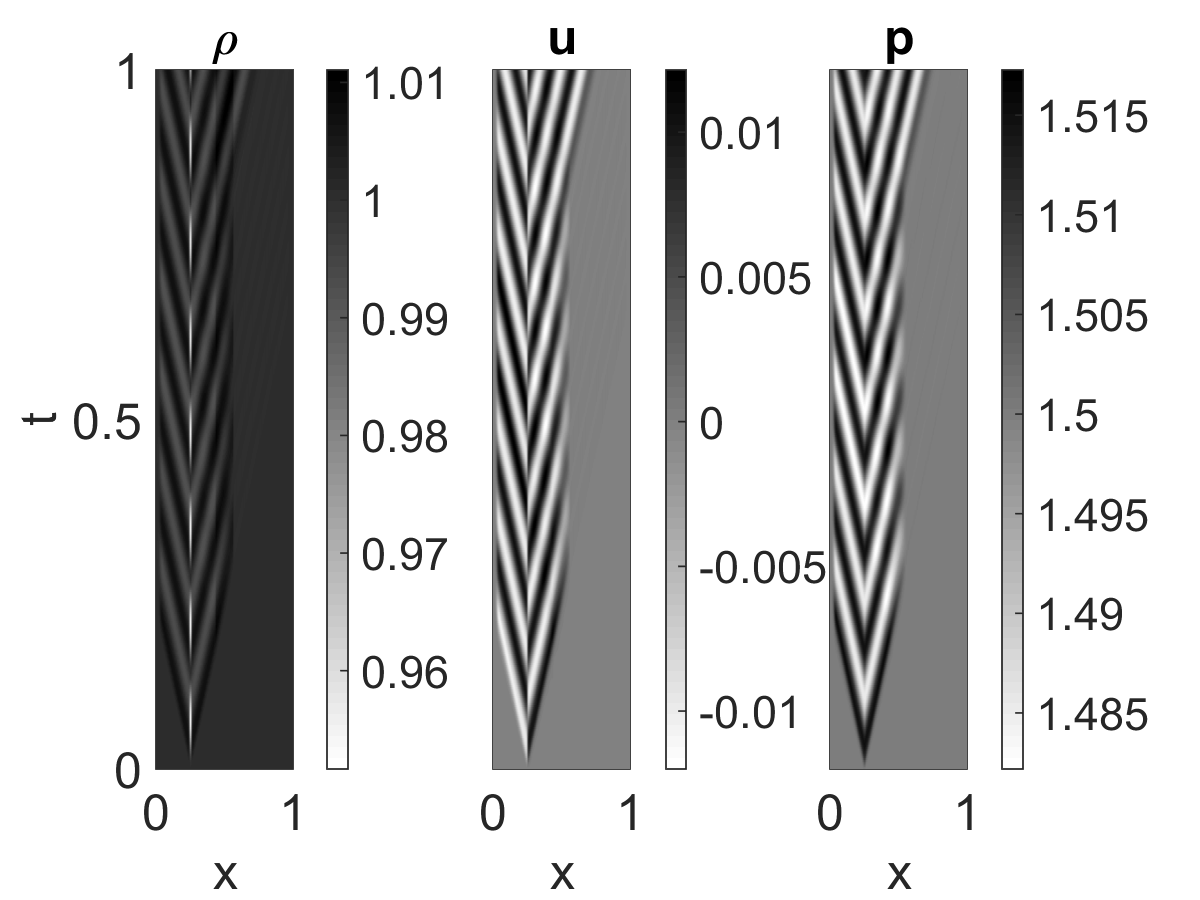}\\
  \includegraphics[width = .49\textwidth]{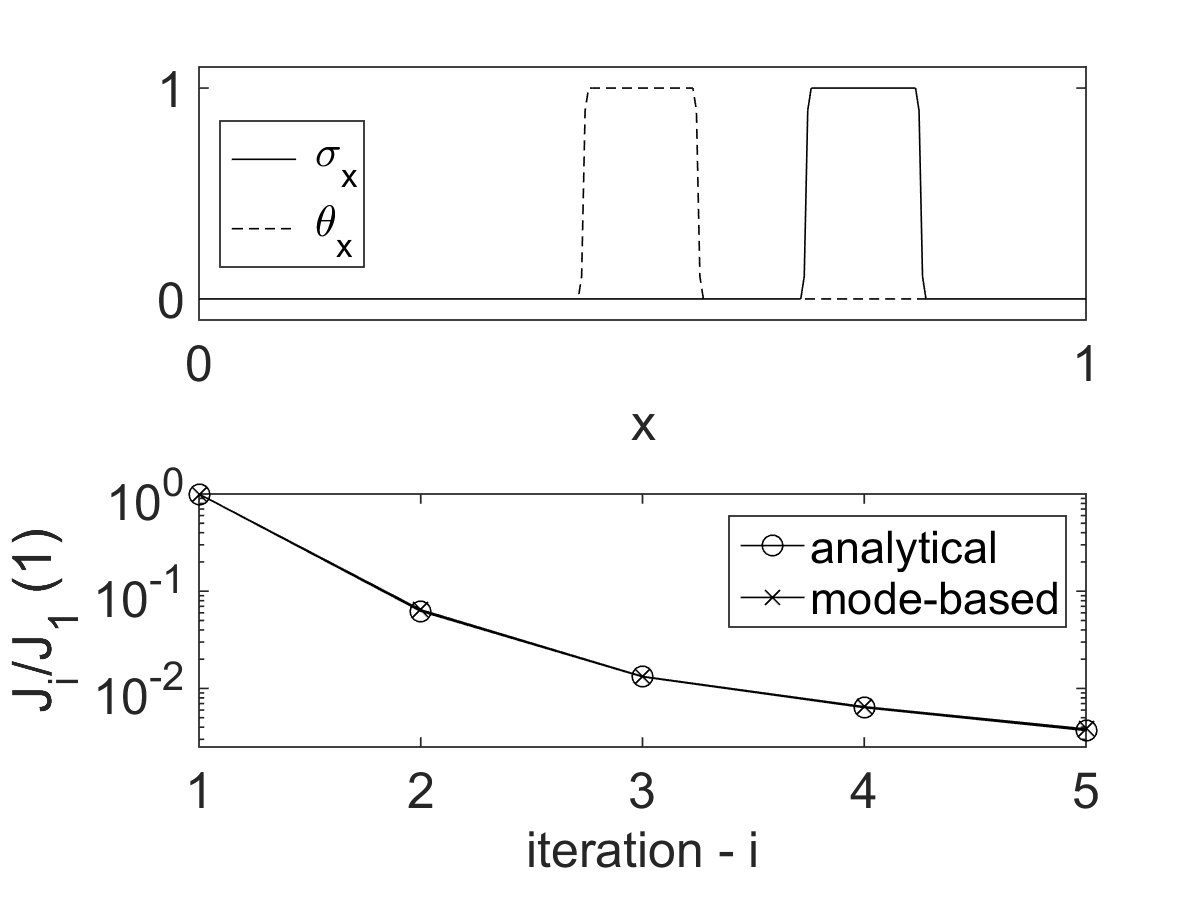}
  \includegraphics[width = .49\textwidth]{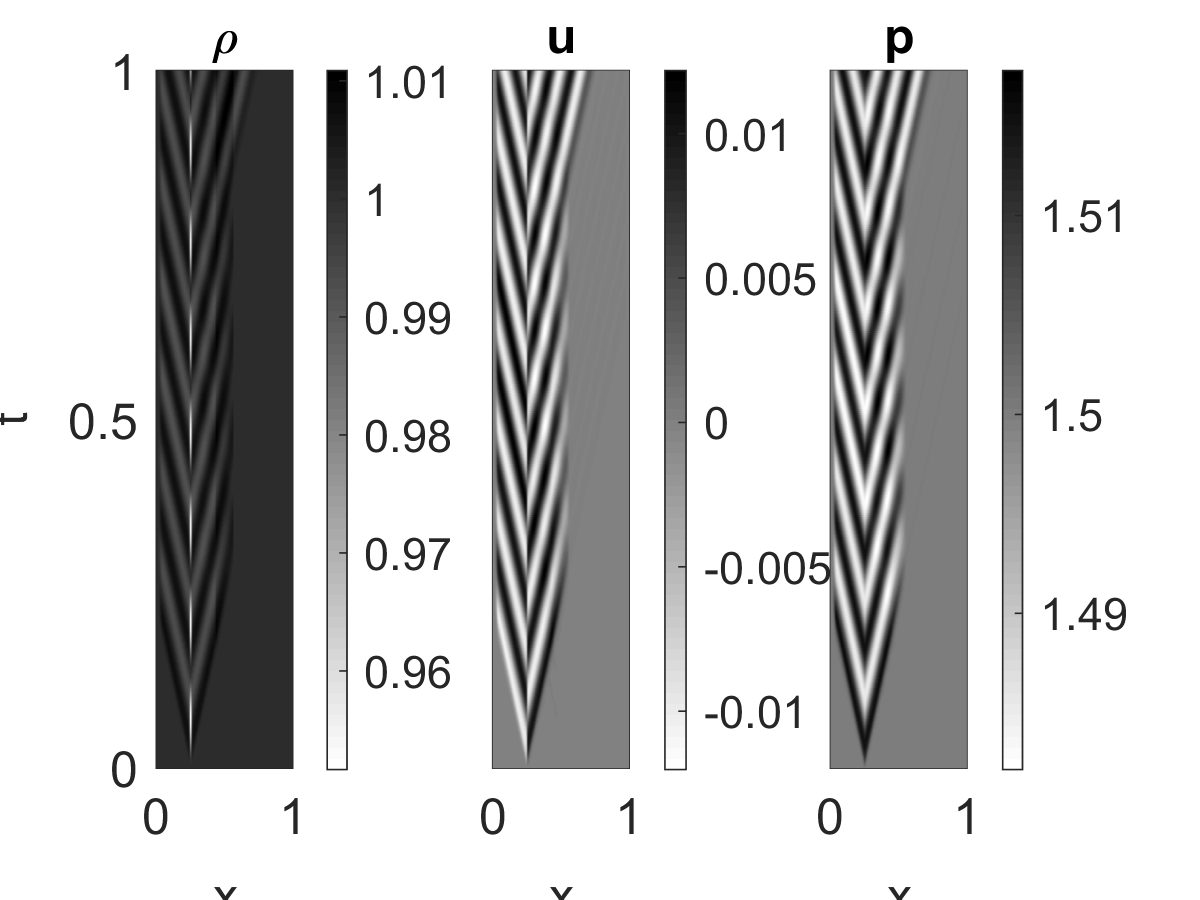}
  \caption{O1 - Noise Cancellation. Solution of the primal system without optimization ($f_p(x,t) = 0$) (top-left) and primal solution after five iteration using the analytical adjoint (top-right).
	  Optimization setup and progress of the objective function normalized with respect to the first iteration (bottom-left) and primal solution after five iteration using the mode-based adjoint  (bottom-right).
  \label{fig_results_c8}
  }
\end{figure}

In order to demonstrate the applicability of the mode-based adjoint approach for optimization tasks, the previously discussed setup (E4) is modified.
The computational domain is extended to a total length of $L = 4\pi$ resolved by 256 equidistantly distributed points.
In summary 384 time steps are simulated at a CFL-condition of $0.75$.
The system is excited by means of a harmonic pressure source with a frequency of $f = 0.75$ Hz located at $x_s = L/4$.
The resulting flow field is characterized by acoustics waves as shown in Fig.~\ref{fig_results_c8} (top-left).

The overall target of the optimization is to minimize the integral objective
\begin{equation}
  J = \iint \left( p(x,t) - p_{\mathrm{target}}\right)^2 \sigma_x ~\mathrm d t
\end{equation}
by means of an adjoint-based adaptation of a source term $f_p$ in the right-hand-side of the pressure equation in \eqref{eq_euler}, with $p_{\mathrm{target}} = 1.5$ as target pressure and $\sigma_x$ a spatial weight, defined by a Gauss-smoothed step function located at $3/4 L$ of the computational domain length, see Fig.~\ref{fig_results_c8}.
According to \eqref{eq_adjoint_sensitivity} the adjoint of \eqref{eq_euler} provides the gradient of $J$ with respect to $f_p$.
The solution is from a human perspective trivial, as in the middle region $\theta_x$ a source is created, which annihilates the left-running acoustic.     
However, several thousand degrees of freedom are adapted. 

The objective is minimized iteratively.
Starting from an initial guess for $f_p(x,t) = 0$ the primal system is solved.
The adjoint system, driven by the term $g$, see \eqref{eq_adjoint_system}, is solved subsequently.  
Based on the solution the force $f_p$ is adapted corresponding to
\begin{align}
  f_p^{n+1} & =  f_p^n + \alpha \left( \frac{\delta J }{\delta f }\right)   \theta_x \no \\ 
 &= f_p^n + \alpha p^* \theta_x 
\end{align}
with a suitable fixed step-width of $\alpha = 2.5$ and the Gaussian-smooth weight $\theta_x$ around $L/2$, which controls the location of the anti-sound-source and reduces to a source in $p$, see Fig.~\ref{fig_results_c8}.
With the updated forcing the primal system is solved again.
The procedure is repeated five times using an analytical adjoint solution and the mode-base approach.

The resulting primal solutions, in which the noise is canceled out by the adapted $f_p$, are shown in Fig.~\ref{fig_results_c8} (right).
There are no identifiable differences between both approaches.
Also the progress of the objective function with respect to the iteration is almost identical.
In both cases the objective is reduced by more than two orders of magnitude.
A nearly invisible deterioration of the convergence using the mode-based approach is negligible.

Please note, that for the dynamic Arnoldi the same calculation plan as in (E4) is employed.
This is particularly remarkable, because, here, for the first time a source term $g$ was present in the adjoint equation and that the plan was not trained to this particular optimization task.

%% file: inctex/conclusion.tex
A new approach for constructing the discrete adjoint operators is presented.
It does not involve the discretization of analytically derived adjoint equations or automatic differentiation.
It is based on the evaluation of the governing primal right-hand-side only. 
Thus, the so-called mode-based adjoint promises to make the adjoint method available for the majority of already existing codes treating various problems with minimal effort.

The approach builds on a modal decomposition similar to an Arnoldi factorization. 
It was found to be a challenge, that in systems of coupled differential equations, like the Euler equations, different variables influence each other in a non-symmetric fashion. 
This is reflected in structural properties of the modes produced by the primal problem,  
being inadequate to describe the adjoint system.      
The dynamic Arnoldi Method (DAM) allowed to overcome this structural problem. 
Thus, the  DAM is an important co--product of this research and might be useful for other problems.

We successfully applied the method to numerous test cases for the Burgers and Euler equations as well as  an optimization problem, showing the principle usability of the method. 
\medskip 

Further investigations are necessary to pave the way towards a  broad applicability. 
This is on one hand an improved training method and on the other hand a generic procedure for boundary treatment.

%% file: inctex/app_calcplans.tex
\clearpage
\twocolumn

\section{Calculation Plans \label{secCalcPlan}}
\begin{table}[h]
  \paragraph{\textsc{Burgers} Equation \hfill \vspace{1em}}
  \begin{tabular}[]{cccc}
  \textbf{Case}				& \textbf{input}	& \textbf{item}	& \textbf{map}	\\
  \hline
  B1 	& I	& 1	& 1 \\
	  & P     & 1     & 1 \\[.5em]

  B2	& I     & 1     & 1 \\
	  & I     & 3     & 1 \\
	  & P     & 2     & 1 \\
	  & P     & 3     & 1 \\[.5em]

  B3	& I     & 1     & 1 \\
	  & I     & 3     & 1 \\
	  & P     & 1     & 1 \\
	  & P     & 3     & 1 \\
	  & P     & 2     & 1 \\
	  & P     & 5     & 1 \\
	  & P     & 6     & 1 \\
  \end{tabular}
  \caption{Calculation plans for the Burgers equation tests using a quality criterion of $10^{-5}$ evaluated over the full time span using each fifth step.\label{app_tbl_calcplans_burgers}}
\end{table}

\vfill\eject
~\vspace{1em}
\begin{table}[h]
  \paragraph{\textsc{Euler} Equations \hfill \vspace{1em}}
  \begin{tabular}[]{cccc}
  \textbf{Case}				& \textbf{input}	& \textbf{item}	& \textbf{map}	\\
  \hline
  E1	& I     & 1     & 1     2     3 \\
	  & P     & 1     & 0     2     0 \\
	  & P     & 1     & 0     0     3 \\
	  & P     & 1     & 0     0     1 \\[.5em]

  E2	& I     & 1     & 1     2     3 \\
	  & I     & 1     & 0     2     0 \\
	  & P     & 1     & 0     0     1 \\
	  & P     & 2     & 0     2     0 \\[.5em]

  E3	& I     & 1     & 1     2     3 \\
	  & P     & 1     & 0     2     0 \\
	  & P     & 1     & 0     0     3 \\
	  & P     & 1     & 0     1     0 \\
	  & P     & 1     & 0     0     1 \\
	  & I     & 1     & 0     0     3 \\[.5em]

  E4/O1	& I     & 1     & 1     2     3 \\
	  & P     & 1     & 0     0     1 \\
	  & P     & 1     & 0     2     0 \\
	  & P     & 1     & 0     0     3 \\
	  & P     & 2     & 0     2     0 \\
	  & I     & 1     & 0     0     3 \\
	  & P     & 3     & 0     2     0 \\
	  & P     & 4     & 0     2     0 \\
	  & P     & 5     & 0     2     0 \\
	  & P     & 6     & 0     2     0 \\
  \end{tabular}
  \caption{Calculation plans for the Euler equations tests using a quality criterion of $10^{-5}$ evaluated over the full time span using each fifth step.\label{app_tbl_calcplans_euler}}
\end{table}